\documentclass[reqno]{amsart}
\usepackage{latexsym}
\usepackage{amssymb, color}
\usepackage{stmaryrd}
\usepackage{epsfig}
\theoremstyle{plain}
\usepackage{txfonts}
\usepackage[sc]{mathpazo}
\def\Xint#1{\mathchoice
   {\XXint\displaystyle\textstyle{#1}}%
   {\XXint\textstyle\scriptstyle{#1}}%
   {\XXint\scriptstyle\scriptscriptstyle{#1}}%
   {\XXint\scriptscriptstyle\scriptscriptstyle{#1}}%
   \!\int}

\def\XXint#1#2#3{{\setbox0=\hbox{$#1{#2#3}{\int}$}
     \vcenter{\hbox{$#2#3$}}\kern-.5\wd0}}

\newcommand{\chara}{1\!\!1}

\newcommand{\cha}{1\!\!1}

\newcommand{\na}{\nabla}

\newcommand{\Pb}{\mathbb{P}}

\newcommand{\ri}{\tau}

\newcommand{\wii}{\ti{w}}

\newcommand{\lt}{\left}
\newcommand{\rt}{\right}
\newcommand{\nl}{\newline}

\newcommand{\nn}{\nonumber}

\newcommand{\lm}{\lambda}

\newcommand{\qd}{\quad}

\newcommand{\ep}{\epsilon}

\newcommand{\wt}{\widetilde}

\newcommand{\N}{\mathrm{I\!N}}
\newcommand{\II}{\mathcal{I}}

\newcommand{\AI}{\mathcal{A}}

\newcommand{\GI}{\mathcal{I}}
\newcommand{\GII}{\mathcal{G}}
\newcommand{\DI}{\mathcal{D}}
\newcommand{\HI}{\mathcal{H}}
\newcommand{\EI}{\mathcal{E}}
\newcommand{\SI}{\mathcal{S}}

\newcommand{\FI}{\mathcal{F}}

\newcommand{\JI}{\mathcal{J}}

\newcommand{\CI}{\mathcal{C}}

\newcommand{\UUI}{\mathcal{U}}
\newcommand{\OI}{\mathcal{O}}
\newcommand{\BI}{\mathcal{U}}

\newcommand{\ti}{\tilde}

\newcommand{\R}{\mathrm {I\!R}}

\newcommand{\ca}[1]{\mathrm{Card}\lt(#1\rt)}
\newcommand{\dia}{\diamondsuit}

\newtheorem{a1}{Lemma}
\newtheorem{a2}{Theorem}
\newtheorem{a3}{Conjecture}

\newtheorem{a5}{Proposition}

\theoremstyle{remark}

\begin{document}
\title[On functions whose symmetric part of gradient agree]
{On functions whose symmetric part of gradient agree and a generalization of 
Reshetnyak's compactness theorem}
\author{Andrew Lorent}
\address{Mathematics Department\\University of Cincinnati\\2600 Clifton Ave.\\ Cincinnati OH 45221 }
\email{lorentaw@uc.edu}
\maketitle

\begin{abstract}
We consider the following question: Given a connected open domain $\Omega\subset \R^n$, suppose 
$u,v:\Omega\rightarrow \R^n$ with $\det(\na u)>0$, $\det(\na v)>0$ a.e.\ are such that 
$\na u^T(x)\na u(x)=\na v(x)^T \na v(x)$ a.e.\ , does this imply a global relation of the 
form $\na v(x)= R\na u(x)$ a.e.\ in $\Omega$ where $R\in SO(n)$? If $u,v$ are $C^1$ it is an exercise to see 
this true, if $u,v\in W^{1,1}$ we show this is false. In Theorem \ref{T1} we prove 
this question has a positive answer if $v\in W^{1,1}$ and $u\in W^{1,n}$ is a mapping of $L^p$ integrable 
dilatation for $p>n-1$. These conditions are sharp in two dimensions and this result represents a 
generalization of the corollary to Liouville's theorem that states that the differential inclusion $\na u\in SO(n)$ can 
only be satisfied by an affine mapping. 

Liouville's corollary for rotations has been generalized by Reshetnyak who proved convergence of 
gradients to a fixed rotation for any weakly converging sequence $v_k\in W^{1,1}$ for which 
$$
\int_{\Omega} \mathrm{dist}(\na v_k,SO(n)) dz\rightarrow 0\text{ as }k\rightarrow \infty.
$$ 
Let $S(\cdot)$ denote the (multiplicative) symmetric part of a matrix. In Theorem \ref{TT2} we prove an analogous result to Theorem 
\ref{T1} for any pair of 
weakly converging sequences $v_k\in W^{1,p}$ and $u_k\in W^{1,\frac{p(n-1)}{p-1}}$ (where $p\in \lt[1,n\rt]$ and the sequence $(u_k)$ has 
its dilatation pointwise bounded above by an $L^r$ integrable function, $r>n-1$) that 
satisfy $\int_{\Omega} \lt|S(\na u_k)-S(\na v_k)\rt|^p dz\rightarrow 0$ as $k\rightarrow \infty$ and 
for which the sign of the $\det(\na v_k)$ tends to $1$ in $L^1$. This result contains Reshetnyak's theorem 
as the special case $(u_k)\equiv Id$, $p=1$.

\end{abstract}

Rigidity of differential inclusions under minimal regularity has been a much studied topic. Probably the best 
known problem of this type is the study of the validity of Liouville's theorem \cite{lou} characterizing 
functions $u$ that satisfy the differential inclusion, 
$$
\na u\in CO_{+}(n):=\lt\{\lm R:\lm>0, R\in SO(n)\rt\}.
$$
Liouville's original theorem was proved for $C^4$ mappings in $\R^3$. This was later generalized by Gehring \cite{geh}, 
Reshetnyak \cite{res}, Bojarski and Iwaniec \cite{bi}, Iwaniec and Martin \cite{iwanm2}, \cite{musvyan}.  In even dimensions the minimal regularity 
for Liouville theorem to hold is $u\in W^{1,\frac{n}{2}}(\Omega)$ (examples show no better result is possible). 
In odd dimensions the optimal regularity is unknown but is conjectured to be $\frac{n}{2}$. 

A corollary to Liouville's theorem is that functions whose gradient belongs to $SO(n)$ are affine.  Note that if  
$u\in W^{1,1}(\Omega)$ and $\na u\in SO(n)$ then $\mathrm{div}(\na u)=\mathrm{div}(\mathrm{cof}(\na u))=0$. 
Thus every co-ordinate function of $u$ weakly satisfies Laplace's equation and hence by Weyl's lemma is 
$C^{\infty}$, thus rigidity for this differential inclusion follows by elementary means. However 
under a much weaker assumptions that $u\in SBV(\Omega)$, $\na u\in SO(n)$ a `piecewise' rigidity 
result has been established in \cite{chamgia}.  Note that the rigidity of the differential inclusion 
$\na u\in SO(n)$ for $\na u\in W^{1,1}$ follows as a highly special case of the following `first guess' conjecture.\nl\nl
\bf `First guess' conjecture. \rm \em 
Let $\Omega\subset \R^n$ be a connected open domain, let $u,v\in W^{1,1}(\Omega)$ and $\det(\na u)>0$, $\det(\na v)>0$ for a.e.\ with 
$$
\na u(x)^T \na u(x)=\na v(x)^T \na v(x)\text{ for }a.e.\ x\in \Omega
$$
then there exists $R\in SO(n)$ such that $\na v=R \na u$ a.e.\rm \nl

As we will show in Example 1, Section \ref{counter}, this conjecture is false. One of the principle aims of this 
paper will be to establish sufficient regularity assumptions required for the 
differential equality $\na u^T \na u=\na v^T \na v$ to imply $\na u=R \na v$ for some $R\in SO(n)$. 
If $u,v\in C^1$ this property would be easy to prove, for $W^{1,1}$ it is not true. In Theorem 
\ref{T1} below we establish the validity of this conjecture with respect to a condition that 
is sharp in two dimensions.
\begin{a2}
\label{T1}
Let $\Omega\subset \R^n$ be a connected open domain, let 
$v\in W^{1,1}(\Omega:\R^n)$ and $u\in W^{1,n}(\Omega:\R^n)$, $\det(\na v)>0$, $\det(\na u)>0$ a.e.\ and $\|\na u(x)\|^n\leq K(x)\det(\na u(x))$ for $K\in L^{P_n}$ where 
\begin{equation}
\label{pnp1}
P_n:=\lt\{\begin{array}{ll} 1&\text{ for } n=2\\
>n-1& \text{ for }n\geq 3\end{array}\rt..
\end{equation}
Suppose 
\begin{equation}
\label{bvbb51}
\na u(x)^T \na u(x)=\na v(x)^T \na v(x) \text{ for }a.e.\ x\in \Omega
\end{equation}
then there exists $R\in SO(n)$ 
\begin{equation}
\label{bvbb50}
\na v(x)=R\na u(x)\text{ for }a.e.\ x\in \Omega.
\end{equation}
\end{a2}

Much interest in the differential inclusion $\na u\in SO(n)$ comes from recent powerful generalization of the corollary to Liouville's theorem that has been established 
in Theorem 3.1 \cite{fmul}. Specifically the $L^2$ distance of the 
gradient of a function away from a \em fixed \em rotation was shown to be bounded by a constant multiple of the $L^2$ distance of the 
gradient away from the \em set \rm of rotations\footnote{A straightforward adaption of 
the proof of Theorem 3.1, \cite{fmul} establishes the same result for $L^p$ control, where $p>1$.}. Previously strong partial results controlling the function 
(rather than the gradient) have been established by John \cite{john}, Kohn \cite{kohn}. 
There has been much work generalizing Theorem 3.1 of \cite{fmul}, for example \cite{chaudmu}, \cite{fac}, \cite{lor5}, \cite{conti1}, \cite{JL}, \cite{contich}.  Part of the motivation for Theorem \ref{T1} is to open a new direction of generalization of Theorem 3.1 of \cite{fmul}. For a 
quantitative version of Theorem \ref{T1} in two dimensions using quite different methods see our companion paper \cite{lor14}. 

Prior to the advances made in \cite{fmul} the most general result generalizing Liouville's theorem 
for mappings with gradient in the space of rotations that gave some control of the gradient was due to Reshetnyak \cite{res}, we state his 
theorem for bounded connected domains.
\begin{a2}[Reshetnyak 1967] 
\label{RTH}

Let $\Omega$ be an open connected and bounded set. If $v_k$ converges weakly in $W^{1,1}(\Omega:\R^n)$ and 
\begin{equation}
\label{dgd1}
\int_{\Omega} \mathrm{dist}(\na v_k,SO(n)) dx\rightarrow 0\text{ as }k\rightarrow \infty
\end{equation}
then $\na v_k$ converges strongly in $L^1$ to a single matrix in $SO(n)$.
\end{a2}

Reshetnyak's Theorem is an example of a result in the more general theory of stability of approximate differential inclusions. Specifically the study of what conditions a set of matrices $K$ must have in order for 
$\int_{\Omega}  \mathrm{dist}(\na v_k,K) dx\rightarrow 0$ to imply $\lt\{\na v_k\rt\}$ is compact in 
$L^1(\Omega)$ for a uniformly bounded Lipschitz sequence, \cite{mulnts}, \cite{tartar1}, \cite{musvyan}. The study of these sets of matrices is closely connected to the theory of quasiconvexity in the calculus of variations \cite{ballq},\cite{morreyq} and was largely motivated by the work of Ball and James \cite{bj}, \cite{bj2}, Chipot and Kinderleher \cite{chk} on variational models of crystal microstructure.   The main result of our paper is a generalization 
of Reshetnyak's theorem, setting $u_k\equiv Id$, $p=1$ in Theorem \ref{TT2} below we recover Theorem 
\ref{RTH}. In the statement of the theorem and from this point on we let $S(\cdot)$ denote the (multiplicative) symmetric part of a matrix. Let $\mathrm{sgn}\lt(\cdot\rt)$ denote the sign of a number, i.e.\ $\mathrm{sgn}(x)=\frac{x}{\lt|x\rt|}$ for $x>0$ and 
$-1$ otherwise.

\begin{a2}
\label{TT2}
Let $\Omega\subset \R^n$ be open and connected and bounded. Let $p\in \lt[1,n\rt]$, $q=\frac{p(n-1)}{p-1}$. Suppose $v_k\in W^{1,p}(\Omega:\R^n)$ converges weakly in $W^{1,p}$ with $\mathrm{sgn}(\det(\na v_k))\overset{L^1}{\rightarrow} 1$ and $u_k\in W^{1,q}(\Omega:\R^n)$ 
converges weakly in $W^{1,q}$, satisfies $\det(\na u_k)>0$ a.e.\ and 
$\|\na u_k\|^n\leq K\det(\na u_k)$ for all $k$ 
where $K\in L^{P_n}$ and $P_n$ satisfies (\ref{pnp1}). If 
\begin{equation}
\label{opps901}
\int_{\Omega} \lt|S(\na v_k)-S(\na u_k)\rt|^p dx \rightarrow 0\text{ as }
k\rightarrow \infty
\end{equation}
then there exists $R\in SO(n)$ such that 
\begin{equation}
\label{opps900}
\lim_{k\rightarrow \infty} \int_{\Omega} \lt|\na v_k-R \na u_{k}\rt| dx=0.
\end{equation}
\end{a2}
For $p=n$, Theorem \ref{TT2} provides a sharp answer to the question, what is the hypothesis necessary such that two weakly 
converging sequences $\na u_k, \na v_k\in W^{1,n}$ with $\int_{\Omega} \lt|S(\na u_k)-S(\na v_k)\rt|^n dz\rightarrow 0$ have the property that there must exists $R\in SO(n)$ so that $\lim_{k\rightarrow \infty} \int_{\Omega} \lt|\na v_k-R \na u_{k}\rt| dx=0$.

By taking $u_k\equiv u$, $v_k\equiv v$ we see Theorem \ref{TT2} generalizes 
Theorem \ref{T1}. Example 1 from Section \ref{counter} shows the necessity (and sharpness in two 
dimensions) of the condition on $(u_k)$. The condition on $(v_k)$ can also easily be seen to 
be necessary, for example by considering $u_k\equiv Id$, $v_k\equiv v$ where $v$ is a 
non affine Lipschitz mapping with its gradient in the set 
$\lt\{\lt(\begin{smallmatrix} 1 & 0 \\ 0 & 1\end{smallmatrix}\rt),\lt(\begin{smallmatrix} 1& 0 \\ 0 & -1 \end{smallmatrix}\rt)\rt\}$.

Another direction of generalization of Theorem \ref{RTH} was proved by M\"{u}ller, Sverak and Yan \cite{musvyan} who generalized Theorem \ref{RTH} for a weakly converging sequence $u_k\in W^{1,\frac{n}{2}}$ where the set of rotations in (\ref{dgd1}) is 
replaced by the set of conformal matrices $CO_{+}(n)$.

As should seem likely from the assumptions of Theorems \ref{TT2} we will be using the powerful 
results established by Iwaniec and Sverak \cite{iws1}, Villamore and Manfredi \cite{manfrediv}, Koskela and Heinonen \cite{kosheim} on 
functions of \em integrable dilatation\rm.  These are functions $u$ for which $L(x):=\frac{\|\na u(x)\|^n}{\det(\na u(x))}$ is a 
positive $L^p$ integrable function, if $L$ is merely positive and finite a.e.\ we say $u$ is a mapping of \em finite dilatation\rm. 
Following \cite{iws1} there has been a well known conjecture that if 
$u$ is a mapping of finite dilatation where $L\in L^{n-1}$ then $u$ is open and discrete. The best 
result to date has been established by  Villamore and Manfredi \cite{manfrediv} whose proved the conjecture 
for functions that satisfy $L\in L^p$ for $p>n-1$. If the conjecture was true for $L\in L^{n-1}$ then Theorem \ref{TT2} would hold for $K\in L^{n-1}$. It is however not clear for $n\geq 3$ if this is the optimal result.

\subsection*{On sharpness} The counter example to the `first guess' conjecture that we construct 
in Section {\ref{counter}} works by squeezing down the center of the square to a point so that 
the interior of the image is disjoint. All known counter examples in higher dimension work in a similar way. If it 
turned out that $L^{p}$ (for $p>n-1$) integrability  
of the dilatation $\frac{\|\na u\|^n}{\det(\na u)}$ was a sharp condition to prevent this, it would 
suggest this condition is sharp for Theorem \ref{T1} and Theorem \ref{TT2}. With this in mind, 
in Section \ref{nonocounter} we consider mappings from the cylinder $B_1(0)\times \lt[0,1\rt]$  
such that $u(B_1(0)\times \lt\{0\rt\})$ consists of a point. If it could be shown such mappings exists with   
$\int_{B_1(0)\times \lt[0,1\rt]} \lt(\frac{\|\na u\|^3}{\det(\na u)}\rt)^p dz<\infty$ for $p<2$ and 
$p\sim 2$ then Theorems \ref{T1}, \ref{TT2} would be sharp. 
However in Proposition \ref{nocounter} it is shown that any radial mapping $u$ of the cylinder that squeezes one 
end to a point but for which each 
co-ordinate function is a product of functions in cylindrical polar co-ordinates that are 
monotonic and convex or concave, then 
$\int_{B_1(0)\times \lt[0,1\rt]} \frac{\|\na u\|^3}{\det(\na u)} dz=\infty$. 
Our guess is that 
Theorem \ref{T1}, Theorem \ref{TT2} are not sharp for $n\geq 3$ and we suspect 
these theorems holds true for functions of  integrable dilatation.

\subsection*{Connections with Stylov decomposition and future directions}

It is worth noting that in two dimensions the validity of `first guess conjecture' is a 
special case of a more general question. 

First some 
background, given $w:\Omega\rightarrow \R^2$, $w(x,y)=(u(x,y),v(x,y))$, for 
$z=x+iy$ let 
$\wii(z)=u(x,y)+i v(x,y)$. Note 
$\frac{\partial \wii}{\partial \overline{z}}(z)=\frac{1}{2}(\frac{\partial}{\partial x}+i\frac{\partial}{\partial y})\wii=\frac{1}{2}(u_x-v_y)+\frac{i}{2}(v_x+u_y)$. And $\frac{\partial \wii}{\partial z}(z)=\frac{1}{2}(\frac{\partial}{\partial x}-i\frac{\partial}{\partial y})\wii=\frac{1}{2}(u_x+v_y)+\frac{i}{2}(v_x-u_y)$. Now 
identifying complex numbers with conformal matrices in the standard way 
$\lt[x+iy\rt]_M=\lt(\begin{smallmatrix} x & -y \\ y & x\end{smallmatrix}\rt)$ we have 
that $\na w(x,y)=\frac{1}{2}\lt[\frac{\partial \wii}{\partial \overline{z}}(z)\rt]_M
\lt(\begin{smallmatrix} 1 & 0 \\ 0 & -1\end{smallmatrix}\rt)+
\frac{1}{2}\lt[\frac{\partial \wii}{\partial z}(z)\rt]_M$. The Beltrami coefficient $\mu(z)$ of $\ti{w}$ is defined by $\frac{\partial \wii}{\partial \overline{z}}(z)=\mu(z)\frac{\partial \wii}{\partial z}(z)$, so $\mu(z)$ relates the conformal part of $\na w$ to the reflection of the anticonformal part of $\na w$. Note that if 
we let $L$ be an affine map with gradient $\lm \lt(\begin{smallmatrix} \cos\theta & -\sin\theta \\ \sin\theta & \cos\theta\end{smallmatrix}\rt)$ then turning 
$L\circ w$ into a complex function we obtain $\lm(\cos\theta+i\sin\theta)\wii$ and 
the Beltrami coefficient of this function is still $\mu(z)$. In other words the 
Beltrami coefficient does not notice changes in gradient made by scaler multiplication or by rotation. It is also not hard to see that if matrices $A$, $B$ have identical Beltrami coefficient then $AB^{-1}\in CO_{+}(n)$ and thus Beltrami coefficient has 
two components and 
`encodes' the geometry of \em how \rm a matrix deforms a ball but does not 
encode any information about the rotation or the size. The 
symmetric part of the gradient has three components and describes both the 
geometry and the size. It should there for not be a surprise that given matrices $A,B\in \R^{2\times 2}$, if $S(A)=S(B)$ then the Beltrami coefficient also agree. There exists a general factorization result known as 
`Stylov' factorization; specifically  for mappings $u_1,u_2$ of finite dilatation and whose 
Beltrami coefficients agree where $u_1$ is a homeomorphism, there exists holomorphic $\phi$ such that 
$u_1=\phi \circ u_2$ (see Theorem 20.4.19 \cite{astala}). If in addition we know that the $S(\na u_1)=S(\na u_2)$ this 
implies $\lt|\na \phi\rt|\equiv 1$ and therefor $\phi$  is a rotation \footnote{Using this and some methods of this paper a short 
proof of Theorem \ref{T1} in two dimensions can be given.}. For higher 
dimensions there is no `Stylov' decomposition and not only are Theorems \ref{T1}, \ref{TT2} 
about non invertible mappings, the methods we use 
to establish them are of very different. 
It is worth noting however that the nature of the factorization is to relate by a 
conformal mapping any two mappings whose gradients pointwise 
deform the ball with the same geometry, ignoring size and rotation. In higher dimensions 
given matrix $A\in \R^{n\times n}$ if we consider $\frac{S(A)}{\lt|S(A)\rt|}$ this 
matrix encodes geometry ignoring size and rotation, so we could consider two functions $u,v$ with the property that $\frac{S(\na u(x))}{\lt|S(\na u(x))\rt|}=\frac{S(\na v(x))}{\lt|S(\na v(x))\rt|}$ for a.e. $x\in \Omega$ and ask if these two functions are related by a 
conformal mapping. We make the following conjecture; 
\begin{a3}
\label{CJ7}
Suppose $\Omega\subset \R^n$ is a bounded open connected domain and $n\geq 3$. Given $u\in W^{1,n}(\Omega)$, $v\in W^{1,1}(\Omega)$ where 
$\det(\na u)>0$, $\det(\na v)>0$ a.e.\ and $u$ satisfies $\int_{\Omega} \lt(\frac{\|\na u\|^n}{\det(\na u)}\rt)^p dz<\infty$ for some 
$p>n-1$ and 
\begin{equation}
\label{band0.1}
\frac{S(\na u(z))}{\lt|S(\na u(z))\rt|}=\frac{S(\na v(z))}{\lt|S(\na v(z))\rt|}\text{ for }a.e.\; z\in \Omega
\end{equation}
then there exists a Mobius transformation $\Phi$ such that $v=\Phi\circ u$. 
\end{a3}
Given that in two dimensions $\frac{S(A)}{\lt|S(A)\rt|}=\frac{S(B)}{\lt|S(B)\rt|}$ is equivalent 
to the Beltrami coefficients of $A$ and $B$ being equal Conjecture \ref{CJ7} would be 
a generalization to `Stylov' factorization to $n\geq 3$, note however Conjecture \ref{CJ7} is 
not true 
in two dimensions (without the assumption of invertibility) as can easily be seen by the complex functions $z^2$, $z^3$. 
One of the main tools we used to prove Theorems \ref{T1}, \ref{TT2} is the quantitative Liouville theorem for rotations of Friesecke, Muller and James \cite{fmul}. In order 
to prove Conjecture \ref{band0.1} what would be required is a quantitative Liouville theorem for conformal matrices. A weakly quantitative result along these lines has been 
proved by Reshetnyak \cite{res1}, and a much stronger quantitative theorem been proved by Faraco and Zhong \cite{fac} for mappings who 
gradient lies in a compact subset of $CO_{+}(n)$ that excludes $0$. Using these theorems and the methods 
of this paper we plan to establish Conjecture \ref{CJ7} in a forth coming work.\nl

\bf Acknowledgments. \rm I would like to thank Stefan M\"{u}ller 
for showing me the connection between the `First guess' conjecture and Liouville's theorem for rotations and 
suggesting this as a topic of study during my stay at the MPI. In addition I am grateful for some very helpful 
initial discussions. I would also like to thank Jon Bevan for pointing out an error in an earlier version of this paper. Finally I would like to thank the referee for many helpful comments and suggestions.

\section{Proof sketch}

\subsection{Sketch of Theorem \ref{T1}} We will begin by sketching the proof in the simplest 
case for smooth globally invertible $u$ and progressively show how the assumptions can be 
weakened till we arrive at hypothesis of Theorem \ref{T1}.

So first we have 
$C^1$ functions $u,v$ where $u$ is globally invertible. Recall for matrix $A\in \R^{n\times n}$ we 
let $S(A)=\sqrt{A^T A}$ be the symmetric part of $A$ and by polar decomposition we have 
$A=R(A)S(A)$ for some $R(A)\in SO(n)$. Form $w(z)=v(u^{-1}(z))$ and 
note that 
\begin{eqnarray}
\label{opz1}
\na w(z)&=&\na v(u^{-1}(x))(\na u(u^{-1}(x)))^{-1}\nn\\
&=&R(\na v(u^{-1}(x)))\lt(R(\na u(u^{-1}(x)))\rt)^{-1}\in SO(n)\nn
\end{eqnarray}
by the Liouville's theorem its clear there exists $R\in SO(n)$ such that 
$\na w(z)=R$ for all $z\in \Omega$. Thus 
\begin{equation}
\label{opx1}
\na v=R \na u\text{ on }\Omega. 
\end{equation}
and result is established. 

Now it can easily be seen that global invertibility is more than we need for this argument 
above to work, if we merely knew that for \em every \rm $x\in \Omega$ there exists $r_x>0$ 
such that $u\lfloor B_{r_x}(x)$ is injective then we could use the same argument to show there 
exists $R_x\in SO(n)$ such that $R_x\na u=\na v$ on $B_{r_x}(x)$. Fix some $x_0$ and 
let 
$$
\UUI:=\lt\{x\in \Omega: R_{x_0}\na u(x)=\na v(x)\rt\}.
$$
For any $x\in \UUI$ we can show $R_x=R_{x_0}$ and thus $\UUI$ is both open and closed. As $\Omega$ 
is connected it is clear that $\UUI=\Omega$. So if we merely have a set $\GI\subset \Omega$ where $\lt|\Omega\backslash \GI\rt|=0$, $\GI$ is connected and $u$ is locally injective on every point $x\in \GI$ then the argument above 
will still carry through. 

Now suppose $v,u\in W^{1,1}$ and $u$ open and discrete then by a theorem of Chernavskii \cite{chern1} we know that the set 
of points on which $u$ fails to be locally injective (the so called `branch set') which we denoted by 
$B_u$, is a set of topological dimension less than $n-2$. Thus by Example VI 11 p93 \cite{wallman} we know that 
$\Omega\backslash B_u$ is connected. However we are blocked from directly carrying out the previous 
argument by the fact that even if we knew $u^{-1}:u(B_{r_x}(x))\rightarrow B_{r_x}(x)$ has Sobolev regularity it 
does not follow that $w=v\circ u^{-1}$ is defined or if it is defined to what extent some kind of chain rule holds for it. Therefor 
more regularity of 
$u$ is required. If $u$ was quasiregular then $u\lfloor B_{r_x}(x)$ is quasiconformal and 
hence $u^{-1}\lfloor u(B_{r_x}(x))$ is quasiconformal and so $w$ would be a well defined 
Sobolev function and the chain rule holds for $v\circ u^{-1}$. Thus we could show $\na w\in SO(n)$ 
on $u(B_{r_x}(x))$ and the argument could be completed to establish $R\na u=\na v$ on $\Omega$.

Now from the other direction let us consider how Theorem \ref{T1} could fail, take the map $P:Q_1(0)\rightarrow \R^2$ defined by $P(x,y)=(x,xy)$ for $x>0$ and $P(x,y)=(x,-xy)$ for $x>0$. So this 
map takes the unit square and squeezes the center down to form a bow tie. If we take another mapping 
$H$ that leaves the left hand side of the bow tie alone and rotates down 
the right hand side. Then comparing $H\circ P$ and $P$ we have that the symmetric part of the 
gradient of both of these functions agree almost everywhere, however we clearly have that there 
is no rotation $R$ such that (\ref{opx1}) holds true. For more details of this mapping see 
Example 1, Section \ref{counter}. It is easy to see that the dilatation $\frac{\|\na P(x,y)\|^2}{\det(\na P(x,y))}\sim x^{-1}$ and so is not integrable. On the other hand in two dimension from the work of Iwaniec, Sverak \cite{iws1} 
we know mappings of integrable dilatation share many of the strong properties of quasiregular mappings. What is not clear for these mappings is if the chain rule holds for the composition 
$v\circ u^{-1}$, we do however at least know Sobolev regularly of $u^{-1}$ by \cite{henclkoma}.

If we have a Lipschitz function $f$ and a function $g\in W^{1,p}$ by considering the difference 
quotients of $f\circ g$ it is easy to see that $f\circ g\in W^{1,p}$. 
This does not mean that the chain rule holds, however in the case where $\det(\na g(x))>0$ 
for a.e.\ $x$ we can apply the general BV chain rule 
of Ambrosio, Dal Maso \cite{ambdal} . Given this is the 
case a natural approach is for us to consider replacing $v$ with a Lipschitz function $\ti{v}$ 
with the property that $\int \lt|\na v-\na \ti{v}\rt|^p dx\approx 0$. Such a function 
can be found by the now standard truncation arguments via maximal functions of 
\cite{zhang}, \cite{acem}. The difficulty of this approach is that the composed function 
$\ti{v}\circ u$ will not necessarily have its gradient in the set of rotations so the best we 
can hope for is an approximate differential inclusion 
\begin{equation}
\label{opx3}
\int_{u(B_{r_x}(x))} d(\na (\ti{v}\circ u^{-1}), SO(n)) dx\approx 0. 
\end{equation}
By use of the previously mentioned quantitative Liouville theorem of Friesecke, M\"{u}ller and James \cite{fmul} we would then be able to conclude that there exists 
$R\in SO(n)$ such that 
$$
\int_{u(B_{r_x}(x))} \lt|\na (\ti{v}\circ u^{-1})-R\rt| dx\approx 0.
$$ 
We have the following estimates 
\begin{eqnarray}
&~&\int_{u(B_{r_x}(x))} d\lt(\na \lt(\ti{v}(u^{-1}(z))\rt),SO(n)\rt) dz\nn\\
&~&\qd\leq \int_{u(B_{r_x}(x))} 
\lt|\lt(\na \ti{v}(u^{-1}(z))-\na v(u^{-1}(z))\rt)\na u(u^{-1}(z))^{-1}\rt| dz\nn\\
&~&\qd\leq
\int_{B_{r_x}(x)} \lt|\lt(\na \ti{v}(y)-\na v(y)\rt)\mathrm{ADJ}\lt(\na u(y)\rt)\rt| dy.
\end{eqnarray}
So in order to control this expression we need the appropriate integrability assumptions on $\na v$, $\na \ti{v}$ and $\na u$. Since 
$v\in W^{1,p}(\Omega)$ so $\|v-\ti{v}\|_{W^{1,p}}\approx 0$ and so by Holder's inequality we have 
$$
\int_{u(B_{r_x}(x))} d\lt(\na \ti{v}(u(z)),SO(2)\rt) dz\leq \lt(\int_{B_{r_x}(x)} \lt|\na \ti{v}-\na v\rt|^p dz \rt)^{\frac{1}{p}}
\lt(\int_{B_{r_x}(x)} \lt|\na u\rt|^{\frac{(n-1)p}{p-1}} dz \rt)^{\frac{p-1}{p}}\approx 0.
$$
So we can apply Friesecke, M\"{u}ller and James \cite{fmul} and conclude that there exists $R_x\in SO(2)$ such that 
$$
\int_{u\lt(B_{r_x}(x)\rt)} \lt|\na (\ti{v}\circ u^{-1})-R_x\rt| dz\approx 0.
$$
Unwrapping this and taking the limit as $\ti{v}\rightarrow v$ we have that $\na v=R_x \na u$ on $B_{r_x}(x)$ and we can complete the 
argument by showing this relation holds globally off the branch set of $u$.

\subsection{Sketch of Theorem \ref{TT2}} The starting point for Theorem \ref{TT2} is Theorem 1.4
of \cite{geriw} that allows us to conclude that letting $u$ denote the weak limit of $u_k$ we have 
$\frac{\|\na u(z)\|^n}{\det(\na u(z))}\leq K(z)$ for a.e.\ $z\in \Omega$. Let $v$ denote the weak limit of $v_k$. Since 
$u\in W^{1,n}(\Omega)$ and $v\in W^{1,1}(\Omega)$ for a.e.\ $x\in \Omega$ both $u,v$ are approximately differentiable, hence from some $r_x>0$ we have that 
$\frac{\lt|u(z)-(u(x)+\na u(x)(z-x))\rt|}{\lt|z-x\rt|}\approx 0$ and 
$\frac{\lt|v(z)-(v(x)+\na v(x)(z-x))\rt|}{\lt|z-x\rt|}\approx 0$ for all $z\in B_{r_x}(x)$. Now 
as $v_k\overset{L^{1}(\Omega)}{\rightarrow} v$ and 
$u_k\overset{L^{1}(\Omega)}{\rightarrow} u$ so for large enough $k$ we have that $v_k$ and 
$u_k$ are very well approximated by the affine maps $W^v_x(z):=v(x)+\na v(x)(z-x)$ and $W^u_x(z):=u(x)+\na u(x)(z-x)$. Now for large enough $k$ we also know that 
$\Xint{-}_{B_{r_x}(x)} \lt|S(\na u_k)-S(\na v_k)\rt| dz\approx 0$ and thus using Lemma \ref{LL5.5} we have that there exists $R_x\in SO(n)$ such that 
\begin{equation}
\label{opx12}
\Xint{-}_{B_{\frac{r_x}{2}}(x)} \lt|\na v_k-R_x \na u_k\rt| dz\approx 0. 
\end{equation}
By Poincare's inequality for 
some affine map $L_x$ with $\na L_x=R_x$ we have $\Xint{-}_{B_{\frac{r_x}{2}}(x)} \lt|v_k-L_x\circ u_k\rt| dz\approx 0$. Recall $v_k$ and $u_k$ are very well approximated by $W^v_x$ and $W^u_k$ thus it must follow that $\na v(x)=R_x \na u(x)$. This implies $S(\na v(x))=S(\na u(x))$ for a.e.\ $x\in \Omega$ and hence we are in a position to apply Theorem \ref{T1}. Thus there exists $R\in SO(n)$ such that 
$\na u(x)=R \na v(x)$ for a.e.\ $x\in \Omega$. Now again as $v_k$ and $u_k$ are $L^{\infty}$ close 
to $v$, $u$ by Poincare's inequality from (\ref{opx12}) we have that $R_x\approx R$. By covering $\Omega$ with a not too overlapping collection $\lt\{B_{r_{x_1}}(x_1), B_{r_{x_2}}(x_2), \dots B_{r_{x_q}}(x_q)\rt\}$ we have that for each $i$, $R_{x_i}\approx R$ and so 
$\int_{\Omega} \lt|\na v_k-R\na u_k\rt| dx\approx 0$ for all large enough $k$. 

Given the similarity between Lemma \ref{LL5.5} and Lemma \ref{LL5} it may seem curious that we need Lemma \ref{LL5.5} at all. The reason 
is that the estimate in Lemma \ref{LL5} gets control of $\lt|\na u_k- r_k \na v_k\rt|$ on a ball of radius 
$c r \exp\lt(-\frac{A^r_k(x)}{\EI^n}\rt)$ where $A_k^r(x)=\Xint{-}_{B_{r}(x)} \lt|\na u_k\rt|^n dx$. In order to obtain global control 
of $\lt|\na u_k- R\na v_k\rt|$ for some fixed $R\in SO(n)$ over the whole of some (large) subset $\Omega'\subset \Omega$ we would 
need a collection $\lt\{B_{c r_q \exp\lt(-\frac{A_k^{r_q}(x_q)}{\EI^n}\rt)}(x_q):x_q\in \Omega'\rt\}$ for which 
\begin{equation}
\label{df1}
\sum_{q} \chara_{B_{c r_q \exp\lt(-\frac{A_k^{r_q}(x_q)}{\EI^n}\rt)}(x_q)}\leq 5. 
\end{equation} 

For this to work, (i.e.\ to be able to apply Lemma \ref{LL5}) we would need an estimate of the form $\sum_{q} A_k^{r_q}(x_q) (r_q)^n\leq c$ which would be available by equi-integrability  
of $\lt|\na u_k\rt|^n$ if $\lt\{B_{r_q}(x_q):q\in \mathbb{N}\rt\}$ did not overlap by some fixed constant. However this completely 
fails to be a consequence of (\ref{df1}) and so no such estimate is available and more subtle arguments are needed to first establish Lemma 
\ref{LL5.5} and get control of the functions in a ball of radius $\frac{r_q}{2}$ and then in the proof of Theorem \ref{TT2} to carefully check the hypothesis of this Lemma \ref{LL5.5} are satisfied.

%
%
%
%

\section{Proof of Theorem \ref{T1}}

\begin{a1}
\label{LL5} Let $p\in \lt[1,n\rt]$, $q=\frac{p(n-1)}{p-1}$. 
Suppose $v\in W^{1,p}(B_r(x):\R^n)$ and $u\in W^{1,q}(B_r(x):\R^n)$ is a homeomorphism of integrable dilatation, i.e. 
there exists positive function $K\in L^1$ such that $\|\na u(z)\|^n\leq K(z)\det(\na u(z))$. Suppose   
for some constant $\EI\in (0,1)$  
\begin{equation}
\label{band040}
B_{\EI r}(u(x))\subset u\lt(B_{\frac{r}{4}}(x)\rt)
\end{equation}
and for $\ep>0$ such that $r^{-32n\lt(n-1\rt)}\leq \mathrm{In}(2+\ep^{-\frac{1}{4}})$ we have 
\begin{equation}
\label{chreqq2}
\Xint{-}_{B_r(x)} \lt|S(\na u)-S(\na v)\rt|^p dz\leq \ep
\end{equation}
and
\begin{equation}
\label{xband0701}
\Xint{-}_{B_r(x)} \lt|\mathrm{sgn}(\det(\na v))-1\rt| dz \leq \ep
\end{equation}
then for positive constants $\CI_0=\CI_0(n)$, $\CI_1=\CI_1(n,\int_{B_r(x)} K dz)$ there exists $R\in SO(n)$ such that 
\begin{equation}
\label{ghf27}
\Xint{-}_{B_{\CI_0 r \exp\lt(-\frac{A_u^r}{\EI^n}\rt)}(x)} 
\lt|\na u-R \na v\rt| dz\leq \CI_1 A^r_u \exp\lt(\frac{n A_u^r}{\EI^n}\rt) \lt(\mathrm{In}(2+\ep^{-\frac{1}{2}})\rt)^{-\frac{1}{32n}},
\end{equation} 
for 
\begin{equation}
\label{ghf324}
A_u^r:=\Xint{-}_{B_r(x)} \lt|\na u\rt|^q dz.
\end{equation}

\end{a1}
\em Proof of Lemma \ref{LL5}. \rm First some notation, given subset 
$S$ of $\R^n$ or $\R^{n\times n}$ and $h>0$ let 
\begin{equation}
\label{nna1}
N_h(S):=\lt\{X:\inf\lt\{\lt|X-Y\rt|:Y\in S\rt\}<h\rt\}.
\end{equation} 


Note $q\geq p$ so by Holder, (\ref{chreqq2}) and (\ref{ghf324}) implies that $\Xint{-}_{B_{\frac{r}{2}}(x)} \lt|\na v\rt|^p dz\leq (A_u^r+1)$. Define 
$$
M_{\gamma}:=\lt\{z\in B_{\frac{r}{4}}(x): \sup_{h\in (0,\frac{r}{4})} \Xint{-}_{B_h(z)} \lt|\na v\rt|^p dz>\gamma  \rt\}.
$$
We have $M_{\gamma_2}\subset M_{\gamma_1}$ for $\gamma_2>\gamma_1$ and $\lt|M_{\gamma}\rt|\rightarrow 0$ as $\gamma\rightarrow 0$. Thus we can find $\lm>0$ large enough so that 
\begin{equation}
\label{ghf29}
\int_{\lt\{z\in B_{\frac{r}{2}}(x): \lt|\na v(z)\rt|>\lm\rt\}} \lt|\na v\rt|^p dz<\sqrt{\ep}r^n 
\end{equation}
and 
\begin{equation}
\label{nmmm9}
\|\na v\|_{L^p(M_{\lm})}\leq c\ep^{\frac{1}{2p}}r^{\frac{n}{p}}.
\end{equation}
Arguing as in Theorem 3, Section 6.6.3 \cite{evans2} we have that 
\begin{equation}
\label{ghf30}
\lt|M_{\lambda}\rt|\leq c\lambda^{-p}\int_{\lt\{z\in B_{\frac{r}{2}}(x):\lt|\na v(z)\rt|>\lm\rt\}} \lt|\na v\rt|^p dz\overset{(\ref{ghf29})}{\leq} c\lambda^{-p}\sqrt{\ep}r^n.
\end{equation} 
Letting $\|\|$ denote the sup norm on the space of matrices,  
$$
\lt|\|\na u\|-\|\na v\|\rt|\leq c\lt|S(\na u)-S(\na v)\rt|
$$ 
we have $\|\|\na u\|-\|\na v\|\|_{L^p(B_r(x))}\overset{(\ref{chreqq2})}{\leq} c \ep^{\frac{1}{p}}r^{\frac{n}{p}}$. So 
\begin{equation}
\label{nmmm20}
\|\na u\|_{L^p(M_{\lm})}\leq c\|\|\na u\|-\|\na v\|\|_{L^p(B_r(x))}+
c\|\na v\|_{L^p(M_{\lm})}\overset{(\ref{nmmm9})}{\leq} c\ep^{\frac{1}{2p}}r^{\frac{n}{p}}.
\end{equation}
By Proposition A1 \cite{fmul} there exists $c\lm$-Lipschitz function $s$ such that 
\begin{equation}
\label{fkf40.5}
\int_{B_{\frac{r}{4}}(x)} \lt|\na v-\na s\rt|^p dz\leq c
\int_{\lt\{z\in B_{\frac{r}{2}}(x): \lt|\na v(z)\rt|>\lm\rt\}} \lt|\na v\rt|^p dz
\end{equation}
And so by (\ref{ghf29}) 
\begin{equation}
\label{npn1}
\int_{B_{\frac{r}{4}}(x)} \lt|\na v-\na s\rt|^p dz\leq c\sqrt{\ep}r^n.
\end{equation}
Let
%
%
\begin{equation}
\label{zband051}
\BI:=\lt\{z\in B_{\frac{r}{4}}(x):v(z)\not=s(z)\rt\}, 
\end{equation}
by Proposition A1 we also have that
\begin{equation}
\label{ghf31}
\BI\subset M_{\lm}.
\end{equation}

\em Step 1. \rm Let $w:B_{\EI r}(u(x))\rightarrow \R^n$ be defined by $w(z):=s(u^{-1}(z))$. There exists $R\in SO(n)$
\begin{equation}
\label{ghf104}
\Xint{-}_{B_{\EI r}(u(x))} \lt|\na w-R\rt| dz
\leq \frac{c \EI^{-n} A^r_u}{\sqrt{\mathrm{In}(2+\ep^{-\frac{1}{2}}})}.
\end{equation}

\em Proof of Step 1. \rm Since $u$ is a mapping of finite dilatation and $\na u\in L^n(\Omega)$ by 
Theorem 1.2. \cite{cohema} we have that $u^{-1}\in W^{1,1}(u(\Omega))$ and $u^{-1}$ is a mapping of finite dilatation. 
Now by the BV chain rule of Ambrosio, DalMaso \cite{ambdal}, (see Theorem 3.101 
\cite{amfupa} or Corollary 3.2 \cite{ambdal}) for a.e. $x\in u(\Omega)$ the restriction of $s$ to the affine space $A(x):=u^{-1}(x)+\lt\{\na u^{-1}(x)v:\in \R^n\rt\}$ is differentiable 
at $u^{-1}(x)$. Since for a.e.\ $x\in u(\Omega)$, $\det(\na u^{-1}(x))>0$ so $A(x)=\R^n$. Thus by Corollary 3.2 \cite{ambdal}, 
$\na w(x)=\na s(u^{-1}(x))\na u^{-1}(x)$. 

Define $\JI$ to be the $n\times n$ diagonal matrix defined by $\JI=\mathrm{diag}(1,1,\dots ,1,-1)$. Let $\II\in \lt\{Id, \JI\rt\}$, note that for 
any $\SI\subset B_{\frac{r}{4}}(x)\backslash \BI$
\begin{eqnarray}
\label{opx400}
&~&\int_{u(\SI)} d(\na w(z),SO(n)\II) dz\nn\\
&~&\qd\qd=\int_{u(\SI)} d\lt(\na s(u^{-1}(z)))\lt(\na u(u^{-1}(z))\rt)^{-1},SO(n)\II\rt) dz\nn\\
&~&\qd\qd\overset{(\ref{zband051})}{=}\int_{\SI} d\lt(\na v(y)\lt(\na u(y)\rt)^{-1},SO(n)\II\rt) \det(\na u(y)) dy.
\end{eqnarray}

Now for any $y\in B_{\frac{r}{2}}(x)$, 
let $R_{v}(y)\in O(n)$, $R_u(y)\in SO(n)$ such that $\na v(y)=R_{v}(y)S(\na v(y))$, $\na u(y)=R_u(y) S(\na u(y))$. Now 
\begin{equation}
\label{nmmm1}
\na v(y)(\na u(y))^{-1}=R_{v}(y) S(\na v(y))(S(\na u(y)))^{-1} R_u(y)^{-1}. 
\end{equation}
So 
\begin{eqnarray}
\label{band05}
\lt|\na v(y)(\na u(y))^{-1}-R_{v}(y) R_u(y)^{-1}\rt|&\leq&c\lt|S(\na v(y))(S(\na u(y)))^{-1}-Id\rt|\nn\\
&=&c\lt|\lt(S(\na v(y))-S(\na u(y))\rt) \lt(S(\na u(y))\rt)^{-1}\rt|\nn\\
&=&c\lt|\lt(S(\na v(y))-S(\na u(y))\rt) \lt(\na u(y)\rt)^{-1}\rt|.
\end{eqnarray}
%
%
Thus 
\begin{eqnarray}
\label{aga}
&~&\int_{B_{\frac{r}{2}}(x)} \lt|\na v(y)(\na u(y))^{-1}-R_{v}(y) R_u(y)^{-1}\rt|\det(\na u(y)) dy\nn\\
&~&\qd\qd\qd\overset{(\ref{band05})}{\leq}c\int_{B_{\frac{r}{2}}(x)} \lt|S(\na v(y))-S(\na u(y))\rt|\lt|\mathrm{ADJ}(\na u(y))\rt| dy\nn\\
&~&\qd\qd\qd \leq c\int_{B_{\frac{r}{2}}(x)} \lt|S(\na v(y))-S(\na u(y))\rt|\lt|\na u(y)\rt|^{n-1} dy\nn\\
&~&\qd\qd\qd\leq c\|S(\na u)-S(\na v)\|_{L^p(B_{\frac{r}{2}}(x))}\lt(\int_{B_{\frac{r}{2}}(x)} \lt|\na u\rt|^q dy\rt)^{\frac{p-1}{p}}\nn\\
&~&\qd\qd\qd\overset{(\ref{chreqq2}),(\ref{ghf324})}{\leq} c\ep^{\frac{1}{p}} A_u^r r^n.
\end{eqnarray}
%
%
%
Let 
\begin{equation}
\label{xband0999}
\DI:=\lt\{z\in B_{\frac{r}{2}}(x):\det(\na v(z))\leq 0\rt\}.
\end{equation}
By (\ref{xband0701})
\begin{equation}
\label{xband0901}
\lt|\DI\rt|\leq \ep r^n.
\end{equation}
Now by (\ref{band05}) and the definition of $\DI$ we know 
$$
d(\na v(y)(\na u(y))^{-1},SO(n))\leq c\lt|(S(\na v(y))-S(\na u(y)))(\na u(y))^{-1}\rt|\text{ for }y\in B_{\frac{r}{2}}(x)\backslash \DI.
$$
Thus taking $\SI=B_{\frac{r}{4}}(x)\backslash (\BI\cup \DI)$ in (\ref{opx400}) for 
the case $\II=Id$ and applying (\ref{aga}) we have 
\begin{equation}
\label{ghf40}
\int_{u(B_{\frac{r}{4}}(x)\backslash (\BI\cup \DI))} d(\na w(z),SO(n)) dz\leq c\ep^{\frac{1}{2n}} A_u^r r^n.
\end{equation}
Now by (\ref{band05})
$$
d(\na v(y)(\na u(y))^{-1}, SO(n) \JI)\leq 
c\lt|(S(\na v(y))-S(\na u(y)))(\na u(y))^{-1}\rt|\text{ for } 
y\in \DI\backslash \BI.
$$
So taking $\II=\JI$, 
$\SI=\DI\backslash \BI$ in (\ref{opx400}) and applying (\ref{aga}) we have 
$$
\int_{u(\DI\backslash \BI)} d(\na w(z),SO(n) \JI) dz\leq c \ep^{\frac{1}{p}} (A_u^r)^{\frac{1}{q}} r^n.
$$
Thus 
\begin{equation}
\label{bband01}
\int_{u(\DI\backslash \BI)} d\lt(\na w(z),SO(n) \rt) dz\leq c\lt|u(\DI\backslash \BI)\rt|+c \ep^{\frac{1}{p}}
(A_u^r)^{\frac{1}{q}} r^n.
\end{equation}
%
%
Let $f(p)=\frac{p(n-1)}{p-1}$, so $f'(p)=-\frac{n-1}{(p-1)^2}<0$ for all $p\in \lt(1,n\rt]$. So 
$f$ is decreasing and $f(n)=n$ so 
\begin{equation}
\label{zveq4}
\frac{p(n-1)}{p-1}>n\text{ for all }p\in (1,n).
\end{equation}
%
%
%
Now by Theorem 1.1. \cite{muller} $\int_{B_{\frac{r}{2}}(x)} \det(\na u(z))\mathrm{In}(2+\det(\na u(z))) dz\leq C_2 A_u^r r^n$ for some constant $C_2=C_2(A,n)$. Let 
$\Upsilon:=\lt\{z\in B_r(x):\det(\na u(z))>\ep^{-\frac{1}{2}}\rt\}$. So 
$\lt|u\lt(\DI\backslash \Upsilon\rt)\rt|\overset{(\ref{xband0901})}{\leq} c \ep^{\frac{1}{2}} r^n$.  Now $\lt|u(\Upsilon)\rt|=\int_{\Upsilon} \det(\na u) dz\leq c A_u^r r^n/\mathrm{In}(2+\ep^{-\frac{1}{2}})$. 
So $\lt|u(\DI)\rt|\leq C_2 A^r_u r^n/\mathrm{In}(2+\ep^{-\frac{1}{2}})$ putting this together 
with (\ref{bband01}) we have $\int_{u(\DI\backslash \BI)} d\lt(\na w(z),SO(n) \rt) dz\leq  c A^r_u r^n/\mathrm{In}(2+\ep^{-\frac{1}{2}})$. So applying this to (\ref{ghf40}) we have 
\begin{equation}
\label{bband02}
\int_{u(B_{\frac{r}{4}}(x)\backslash \BI)} d(\na w(z),SO(n)) dz\leq c A_u^r r^n/\mathrm{In}(2+\ep^{-\frac{1}{2}}).
\end{equation}

%
%
%
Now as 
\begin{equation}
\label{zveq1}
\lt|\BI\rt|\overset{(\ref{ghf31})}{\leq}\lt|M_{\lambda}\rt|\overset{(\ref{ghf30})}{\leq}c\lm^{-p}\sqrt{\ep}r^n 
\end{equation}
so (recalling $w=s\circ u^{-1}$ and $s$ is $c\lambda$-Lipschitz)
\begin{eqnarray}
\label{kjk1}
\int_{u(\BI)} \lt|\na w\rt| dz&\leq& c\lm \int_{u(\BI)} \lt|\lt(\na u(u^{-1}(z))\rt)^{-1}\rt| dz\nn\\
&=&c\lm \int_{u(\BI)} \lt|\lt(\na u(u^{-1}(z))\rt)^{-1}\rt| \det(\na u^{-1}(z))\det(\na u(u^{-1}(z)))dz\nn\\
&=&c\lm \int_{\BI} \lt|\lt(\na u(y)\rt)^{-1}\rt|\det(\na u(y)) dy
\leq c\lm \int_{\BI} \lt|\mathrm{ADJ}(\na u(z))\rt| dz\nn\\
&\leq& c\lm \int_{U} \lt|\na u(z)\rt|^{n-1} dz\leq c\lm \lt(\int_{U} \lt|\na u(z)\rt|^{q}\rt)^{\frac{p-1}{p}}\lt|\BI\rt|^{\frac{1}{p}}\nn\\
&\overset{(\ref{zveq1}),(\ref{ghf324})}{\leq}& c\lm \lt(A_u^r r^n \rt)^{\frac{p-1}{p}} \lt(c\lm^{-p}\sqrt{\ep}r^n\rt)^{\frac{1}{p}}\leq 
c A_u^r \ep^{\frac{1}{2p}}r^n.
\end{eqnarray}
Now $\int_{B_r(x)} \lt|\na u\rt|^q dz\overset{(\ref{ghf324})}{\leq} A_u^r r^n$ where $q=\frac{p(n-1)}{p-1}$. So by (\ref{zveq4}) and Holder's inequality we know that 
\begin{equation}
\label{coreq1}
\int_{B_r(x)} \lt|\na u\rt|^n dz\leq A_u^r r^n. 
\end{equation}
Let $\theta\in (0,1)$ such that 
$\frac{1}{n}=\frac{\theta}{p}+\frac{1-\theta}{q}$. So by the $L^p$ interpolation inequality (see Appendix B2 \cite{evanspde}) 
\begin{eqnarray}
\label{zreq6}
\|\na u\|_{L^n(M_{\lm})}&\leq& \|\na u\|_{L^p(M_{\lm})}^{\theta}\|\na u\|_{L^q(M_{\lm})}^{1-\theta}\nn\\
&\overset{(\ref{ghf324}),(\ref{nmmm20})}{\leq}&\lt(\ep^{\frac{1}{2p}} r^{\frac{n}{p}}\rt)^{\theta} 
\lt(\lt(A_u^r\rt)^{\frac{1}{q}} r^{\frac{n}{q}}\rt)^{1-\theta}\nn\\
&\leq& \ep^{\frac{\theta}{2p}} A_u^r r^{n\lt(\frac{\theta}{p}+\frac{1-\theta}{q}\rt)}\nn\\
&\leq& \ep^{\frac{\theta}{2}}A_u^r r.
\end{eqnarray}
Thus
%
%
\begin{eqnarray}
\int_{u(B_{\frac{r}{4}}(x))} d(\na w(z),SO(n)) dz&=&\int_{u(B_{\frac{r}{4}}(x)\backslash \BI)} d(\na w(z),SO(n)) dz\nn\\
&~&\qd\qd\qd+
\int_{u(\BI)} (\lt|\na w(z)\rt|+c) dz\nn\\
&\overset{(\ref{kjk1}),(\ref{bband02})}{\leq}& c\lt|u(\BI)\rt|
+2C_2 A^r_u r^n/\mathrm{In}(2+\ep^{-\frac{1}{2}}) \nn\\
&\overset{(\ref{ghf31})}{\leq}& c\int_{M_{\lm}} \det(\na u) dz+ 3C_2 A_u^r r^n/\mathrm{In}(2+\ep^{-\frac{1}{2}})  \nn\\
&\overset{(\ref{nmmm20})}{\leq}& 3C_2 A_u^r r^n/\mathrm{In}(2+\ep^{-\frac{1}{2}}). 
\end{eqnarray}

So in particular using (\ref{band040}) 
we have $\Xint{-}_{B_{\EI r}(u(x))} d(\na w(z),SO(n)) dz\leq \frac{3C_2 A_u^r\EI^{-n}}{\mathrm{In}(2+\ep^{-\frac{1}{2}})}$. 
Thus by Proposition 2.6 \cite{conti1} we have that 
$$
\Xint{-}_{B_{\EI r}(u(x))} \lt|\na w-R\rt| dz\leq 
c\mathrm{In}\lt(\frac{\mathrm{In}(2+\ep^{-\frac{1}{2}})}{3C_2 A\EI^{-n}}\rt)\frac{3C_2 A_u^r\EI^{-n}}{\mathrm{In}(2+\ep^{-\frac{1}{2}})}
\leq \frac{c A_u^r\EI^{-n}}{\sqrt{\mathrm{In}(2+\ep^{-\frac{1}{2}})}}. 
$$
This completes the proof of Step 1. \nl

\em Step 2. \rm We will show 
\begin{equation}
\label{pox1}
B_{c r \exp (-\frac{A_u^r}{\EI^n})}(x)\subset u^{-1}\lt(B_{\EI r}(u(x))\rt)
\end{equation}

\em Proof of Step 2. \rm Note by equation (2.5) of the proof of Theorem 1 of \cite{manfredi} we know that for 
any $y\in B_{\frac{r}{2}}(x)$, $h\in \lt(0,\frac{r}{2}\rt]$
\begin{eqnarray}
\label{band045}
\mathrm{osc}_{B_h(y)}u&\leq& c\lt(\log\lt(\frac{r}{2h}\rt)\rt)^{-\frac{1}{n}}
\lt(\int_{B_{\frac{r}{2}}(y)} \lt|\na u\rt|^n dz\rt)^{\frac{1}{n}}\nn\\
&\overset{(\ref{coreq1})}{\leq}& c (A_u^r)^{\frac{1}{n}}r \lt(\log\lt(\frac{r}{2h}\rt)\rt)^{-\frac{1}{n}}.
\end{eqnarray}
We claim 
\begin{equation}
\mathrm{dist}\lt(x,u^{-1}\lt(\partial B_{\EI r}(u(x))\rt)\rt)\geq c r \exp\lt(-\frac{A_u^r}{\EI^n}\rt).
\end{equation}
So see this pick 
$z\in u^{-1}\lt(\partial B_{\EI r}(u(x))\rt)\cap B_{\frac{r}{2}}(x)$, since $u$ is a homeomorphism 
$\EI r=\lt|u(z)-u(x)\rt|
\overset{(\ref{band045})}{\leq} c (A_u^r)^{\frac{1}{n}}r \lt(\log\lt(\frac{r}{2\lt|z-x\rt|}\rt)\rt)^{-\frac{1}{n}}$. 
Thus 
$$
\lt(\log\lt(\frac{r}{2\lt|z-x\rt|}\rt)\rt)^{\frac{1}{n}}\EI \leq c (A_u^r)^{\frac{1}{n}}
$$ 
and so  $\log\lt(\frac{r}{2\lt|z-x\rt|}\rt)\EI^n \leq c A_u^r$ and hence 
$\frac{r}{2\lt|z-x\rt|}\leq c \exp(\frac{A_u^r}{\EI^n})$ and 
finally $cr \exp(-\frac{A_u^r}{\EI^n})\leq  \lt|z-x\rt|$ which 
establishes (\ref{pox1}). \nl

\em Proof of Lemma completed. \rm Note that if matrix $B$ satisfies $\|B\|^n\leq Q\det(B)$ then as 
$\Lambda(B):=\inf\lt\{\lt|Bv\rt|:v\in S^{n-1}\rt\}$. If we let $B_I$ be the smallest number such that 
$\det(B)\leq B_I\Lambda(B)^n$ it is well known (see for example \cite{vaisala2} p44) $B_I\leq Q^{n-1}$ so 
\begin{equation}
\label{fqq1}
\Lambda(B)\geq \frac{\det(B)^{\frac{1}{n}}}{Q^{\frac{n-1}{n}}}.
\end{equation} 
So it is an 
exercise to see 
\begin{equation}
\label{lband02}
\lt|AB\rt|\geq \frac{\lt(\det(B)\rt)^{\frac{1}{n}}\lt|A\rt|}{Q^{\frac{n-1}{n}}n}\geq \frac{Q^{-1}}{n^2}\lt|B\rt|\lt|A\rt|\text{ for any } A\in \R^{n\times n}.
\end{equation}
Recall $u$ is of integrable dilatation and so we have function $K$ such that 
$\|\na u(z)\|^n\leq K(z)\det(\na u(z))$. Let  
\begin{equation}
\label{lband01}
\Pb:=\lt\{z\in B_r(x):\lt|K(z)\rt|\geq \lt(\mathrm{In}(2+\ep^{-\frac{1}{2}})\rt)^{\frac{1}{8(n-1)}}\rt\}.
\end{equation}

Now as $\na w(z)=\na s(u^{-1}(z))(\na u(u^{-1}(z)))^{-1}$. Thus by (\ref{ghf104}) Step 1 
\begin{eqnarray}
\label{lband031}
\frac{c A_u^r r^n}{\sqrt{\mathrm{In}(2+\ep^{-\frac{1}{2}})}}&\overset{(\ref{ghf104})}{\geq}& \int_{B_{\EI r}(x)} \lt|\na w-R\rt| \det(\na u^{-1}(z)) 
\det(\na u(u^{-1}(z))) dz\nn\\
&=&\int_{u^{-1}\lt(B_{\EI r}(x)\rt)} \lt|\na s(z)(\na u(z))^{-1}-R\rt|\det(\na u(z)) dz\nn\\
&\geq &\int_{u^{-1}\lt(B_{\EI r}(x)\rt) \backslash \Pb} 
\lt|\lt(\na s(z)-R\na u(z)\rt)(\na u(z))^{-1} \rt|\det(\na u(z)) dz\nn\\
&\geq&c\int_{u^{-1}\lt(B_{\EI r}(x)\rt)\backslash \Pb} 
\lt|\lt(\na s(z)-R \na u(z)\rt) \mathrm{ADJ}(\na u(z))\rt| dz.
\end{eqnarray}
So by using (\ref{fqq1}) $\|\na u(z)^{-1}\|^n\leq K(z)^{n-1}\det((\na u(z))^{-1})$, so 
$$
\|\mathrm{ADJ}(\na u(z))\|^n\leq K(z)^{n-1}\det(\mathrm{ADJ}(\na u(z))). 
$$
Hence if $z\not \in \Pb$ by (\ref{lband02}), 
\begin{equation}
\label{lband030}
\lt|\lt(\na s(z)-R \na u(z)\rt) \mathrm{ADJ}(\na u(z))\rt|\geq 
n^{-2}\lt(\mathrm{In}(2+\ep^{-\frac{1}{2}})\rt)^{-\frac{1}{8}}\lt|\lt(\na s(z)-R \na u(z)\rt) \rt|\lt|\mathrm{ADJ}(\na u(z))\rt|.
\end{equation}
Thus by (\ref{lband031}), (\ref{lband030})
\begin{equation}
\label{lband03}
\frac{c A_u^r r^n}{\lt(\mathrm{In}(2+\ep^{-\frac{1}{2}})\rt)^{\frac{3}{8}}}\geq \int_{u^{-1}\lt(B_{\EI r}(x)\rt)\backslash \Pb} 
\lt|\na s(z)-R \na u(z)\rt| \lt|\mathrm{ADJ}(\na u(z))\rt| dz.
\end{equation}
Now let $\FI:=\lt\{z\in B_r(x):\lt|\mathrm{ADJ}(\na u(z))\rt|<\lt(\mathrm{In}(2+\ep^{-\frac{1}{2}})\rt)^{-\frac{1}{4}} \rt\}$ so 
\begin{equation}
\label{opxx2}
\int_{u^{-1}\lt(B_{\EI r}(x)\rt)\backslash (\FI\cup \Pb)} 
\lt|\na s(z)-R\na u(z)\rt| dz\leq \frac{c A_u^r r^n}
{\lt(\mathrm{In}(2+\ep^{-\frac{1}{2}})\rt)^{\frac{1}{8}}}.
\end{equation}
For any matrix $A\in \R^{n\times n}$ let $M_{ij}(A)$ is the $i,j$ minor of $A$. Thus 
$\lt|M_{ij}(\na u(z))\rt|<\lt(\mathrm{In}(2+\ep^{-\frac{1}{2}})\rt)^{-\frac{1}{4}}$ for any $z\in \FI$. So 
$\lt|\det(\na u(z))\rt|\leq c\lt(\mathrm{In}(2+\ep^{-\frac{1}{2}})\rt)^{-\frac{1}{4}}\lt|\na u(z)\rt|$. Thus 
\begin{eqnarray}
\label{band051}
\int_{\FI} \lt|\det(\na u(z))\rt| dz&\leq&c\lt(\mathrm{In}(2+\ep^{-\frac{1}{2}})\rt)^{-\frac{1}{4}}\int_{\FI} \lt|\na u(z)\rt| dz\nn\\
&\overset{(\ref{coreq1})}{\leq}&c\lt(\mathrm{In}(2+\ep^{-\frac{1}{2}})\rt)^{-\frac{1}{4}} (A_u^r)^{\frac{1}{n}}r^n.
\end{eqnarray}
Now as $\|\na u(z)\|^n\leq K(z)\det(\na u(z))$ for a.e.\ $z\in B_r(x)$. 
Note by (\ref{lband01}) $\lt|\Pb\rt|\lt(\mathrm{In}(2+\ep^{-\frac{1}{2}})\rt)^{\frac{1}{8(n-1)}}\leq \int K dz\leq c$ and thus 
\begin{equation}
\label{eband01.11}
\lt|\Pb\rt|<c\lt(\mathrm{In}(2+\ep^{-\frac{1}{2}})\rt)^{-\frac{1}{8(n-1)}}.
\end{equation} 
So if $z\not\in \Pb$ then $\|\na u(z)\|^n\leq c\lt(\mathrm{In}(2+\ep^{-\frac{1}{2}})\rt)^{\frac{1}{8(n-1)}}\det(\na u(z))$. Thus  
\begin{eqnarray}
\label{nppn1.7}
\int_{\FI\backslash \Pb}\lt|\na u\rt|^n dz&\leq&c \lt(\mathrm{In}(2+\ep^{-\frac{1}{2}})\rt)^{\frac{1}{8(n-1)}}  \int_{\FI\backslash \Pb} \det(\na u(z)) dz\nn\\
&\overset{(\ref{band051})}{\leq}& c\lt(\mathrm{In}(2+\ep^{-\frac{1}{2}})\rt)^{-\frac{1}{8}} (A_u^r)^{\frac{1}{n}} r^n.
\end{eqnarray}
Now 
\begin{equation}
\label{zqz1}
\int_{\Pb} \lt|\na u\rt| dz\leq c\lt|\Pb\rt|^{\frac{n-1}{n}}\lt(\int_{B_r(x)} \lt|\na u\rt|^n dz \rt)^{\frac{1}{n}}\overset{(\ref{coreq1}),(\ref{eband01.11})}{\leq} c\lt(\mathrm{In}(2+\ep^{-\frac{1}{2}})\rt)^{-\frac{1}{8n}}(A_u^r)^{\frac{1}{n}}r.
\end{equation}
So 
\begin{equation}
\label{zqz2.8}
\int_{\FI\cup \Pb} \lt|\na u\rt| dz\leq \int_{\Pb} \lt|\na u\rt| dz+\int_{\FI\backslash \Pb} \lt|\na u\rt| dz
\overset{(\ref{zqz1}),(\ref{nppn1.7})}{\leq} c\lt(\mathrm{In}(2+\ep^{-\frac{1}{2}})\rt)^{-\frac{1}{16n}} (A_u^r)^{\frac{1}{n^2}}r.
\end{equation}
%
%
%
%
By using Holder's inequality we see 
\begin{eqnarray}
\label{nppn2}
\|\na s\|_{L^1(\FI\cup \Pb)}&\leq& c\|S(\na u)-S(\na v)\|_{L^p(B_{\frac{r}{2}}(x))}+c\|\na u\|_{L^1(\FI\cup \Pb)}
+c\|\na v-\na s\|_{L^p(B_{\frac{r}{4}}(x))}\nn\\
&\overset{(\ref{zqz2.8}),(\ref{npn1}),(\ref{chreqq2})}{\leq}& c\lt(\mathrm{In}(2+\ep^{-\frac{1}{2}})\rt)^{-\frac{1}{16n}} (A_u^r)^{\frac{1}{n^2}} 
r.
\end{eqnarray} 
Thus $\int_{\FI\cup \Pb} \lt|\na s-R\na u\rt| dz\overset{(\ref{zqz2.8}),(\ref{nppn2})}{\leq} 
c\lt(\mathrm{In}(2+\ep^{-\frac{1}{2}})\rt)^{-\frac{1}{16n}}(A_u^r)^{\frac{1}{n^2}} r$. Hence putting this together with (\ref{opxx2}), (\ref{npn1}) we have
%
%
%
%
\begin{equation}
\label{ghf33}
\int_{u^{-1}\lt(B_{\EI r}(x)\rt)} 
\lt|\na v-R \na u\rt| dz\leq  c  A_u^r \lt(\mathrm{In}(2+\ep^{-\frac{1}{2}})\rt)^{-\frac{1}{16n}}r.
\end{equation}
And putting this together with (\ref{pox1}) we have 
\begin{equation}
\label{eqeq30}
\Xint{-}_{B_{cr\exp(-\frac{A_u^r}{\EI^n})}(x)} \lt|\na v-R\na u\rt| dz\leq c r^{-n+1}
\exp\lt(\frac{n A_u^r}{\EI^n}\rt) A_u^r \lt(\mathrm{In}(2+\ep^{-\frac{1}{2}})\rt)^{-\frac{1}{16n}}.
\end{equation}
Note that since $r^{-32n\lt(n-1\rt)}\leq \mathrm{In}\lt(2+\ep^{-\frac{1}{2}}\rt)$ so 
$\lt(\mathrm{In}(2+\ep^{-\frac{1}{2}})\rt)^{-\frac{1}{32n}}\leq r^{n-1}$ so 
putting this together with (\ref{eqeq30}) we have established (\ref{ghf27}).  \nl\nl

%
%

\subsection{Proof of Theorem \ref{T1} completed}
By Theorem 1 \cite{manfrediv} $u$ is a discrete open mapping. Let $B_u$ denote 
the set of points $z\in B_r(x)$ such that $u$ is not locally invertible in any neighborhood containing 
$z$. By definition this is a closed set. 

\em Step 1. \rm We will show $u(B_r(x))\backslash u(B_u)$ is connected. 

\em Proof of Step 1. \rm By a theorem of Chernavskii $B_u$ \cite{chern1}, (also see \cite{vaisala}) $B_u$ has topological dimension at most 
$n-2$. By Example VI 11 p93 \cite{wallman} we know $B_r(x)$ can not be separated by  $B_u$ and 
so $B_r(x)\backslash B_u$ is connected. 

\em Proof of Theorem \ref{T1}. \rm For any $z\in B_r(x)\backslash B_u$ by definition of $B_u$ there exists 
$s_z>0$ such that $u_{\lfloor B_{s_z}(z)}$ is injective, therefor by applying Lemma \ref{LL5} we 
know that for some $r_z\in (0,s_z)$ there exists $R_z\in SO(n)$ such that $\na u=R_z \na v$ on $B_{r_z}(z)$. 

Pick $z_0\in B_r(x)\backslash B_u$. Let $z_1\in B_r(x)\backslash B_u$, $z_1\not =z_0$. 
Since $B_r(x)\backslash B_u$ is connected and open and is therefor arcwise connected so there 
exists a homomorphism $\psi:\lt[0,1\rt]\rightarrow B_r(x)\backslash B_u$ with 
$\psi(0)=z_0$, $\psi(1)=z_1$. For each $z\in B_r(z)\backslash B_u$ let 
$$
\alpha_z:=\sup\lt\{\alpha>0:u_{\lfloor B_{\alpha}(z)}\text{ is injective }\rt\}.
$$
It is clear $\beta=\inf\lt\{\alpha_z:z\in \psi(\lt[0,1\rt])\rt\}>0$ since otherwise 
by compactness $B_u\cap \psi(\lt[0,1\rt])\not =\emptyset$. Let 
\begin{equation}
\label{pkl1}
\GII:=\lt\{h\in \lt[0,1\rt]:\na u(z)=R_{z_0} \na v(z)\text{ for }a.e.\ z\in \bigcup_{\gamma\in \lt[0,h\rt]} B_{\frac{\beta}{2}}(\psi(\gamma))\rt\}.
\end{equation}
It is clear $\GII$ is a closed set, it is also straightforward to see it is open because if 
$h\in \GII$ there exits $\wt{R}\in SO(n)$ such that $\na u(x)=\wt{R}\na v(x)$ for a.e.\ $x\in 
B_{\beta}(\psi(h))$. Since we also know $\na u(z)=R_{z_0} \na v(z)$ for a.e.\ 
$z\in B_{\frac{\beta}{2}}(\psi(h))$ it is clear $\wt{R}=R_{z_0}$ and thus there exists $\delta>0$ 
with $(h-\delta,h+\delta)\subset \GII$. As $\GII$ is open and closed in $\lt[0,1\rt]$ and 
as it is non empty we have 
that $\GII=\lt[0,1\rt]$. In particular this implies that for every $z\in B_r(x)\backslash B_u$,  
$\na u(y)=R_{z_0} \na v(y)$ for a.e.\ $y\in B_{\frac{\beta}{2}}(z)$. Thus 
$\na u(z)=R_{z_0} \na v(z)$ for 
a.e.\ $z\in B_r(x)\backslash B_u$. Since $B_u$ has dimension at most $n-2$ we know $\lt|B_u\rt|=0$ 
there for (\ref{bvbb50}) follows immediately. $\Box$

%
%

\section{Preliminary lemmas for Theorem \ref{TT2}}

%
%

\begin{a1}
\label{LL5.5}
Let $r\in (0,1)$, $A>1$, $p\in \lt[1,n\rt]$, $q=\frac{p(n-1)}{p-1}$, $\EI\in (0,1)$. Suppose $v\in W^{1,p}(B_r(x):\R^n)$ and $u\in W^{1,q}(B_r(x):\R^n)$ is a homeomorphism of integrable dilatation. There exists small constant 
$\ep_0=\ep_0(A,r,\EI)$ such that if functions $u,v$ satisfy 
\begin{equation}
\label{zband03}
\Xint{-}_{B_r(x)} \lt|S(\na u)-S(\na v)\rt|^p dz\leq \ep,
\end{equation}
\begin{equation}
\label{zband03.5}
\Xint{-}_{B_r(x)} \lt|\mathrm{sgn}(\det(\na v(z)))-1\rt| dz\leq \ep
\end{equation}
\begin{equation}
\label{zband02}
\Xint{-}_{B_r(x)} \lt|\na u\rt|^q dz\leq A
\end{equation}
for $\ep<\ep_0$ and there exists $\Xi\subset B_{\frac{r}{2}}(x)$ such that 
\begin{equation}
\label{eqeq100}
\lt|B_{\frac{r}{2}}(x)\backslash \Xi\rt|\leq \frac{\CI_0^n}{16^n} r^n \exp\lt(-\frac{n A}{\EI^n}\rt)
\end{equation} 
and 
\begin{equation}
\label{vvv2}
B_{\EI h}(u(x))\subset u(B_{\frac{h}{4}}(x))\text{ for any }x\in \Xi, h\in \lt[\frac{\CI_0}{8} r\exp\lt(-\frac{A}{\EI^n}\rt) ,\frac{r}{2}\rt].
\end{equation}
Then there exists $\CI_2=\CI_2(n,\int_{B_r(x)} K dz)$ and $R\in SO(n)$ such that 
$$
\int_{B_{\frac{r}{2}}(x)} \lt|\na u- R\na v\rt| dz\leq \CI_2 \EI^{-n} 
A^{3n} \exp\lt(\frac{2^{n+2} n^3 A}{\EI^n}\rt) \lt(\mathrm{In}\lt(2+\frac{\ep^{-\frac{1}{2}}}{2^{n}}\rt)\rt)^{-\frac{1}{64 n^2}} r^{n}.
$$
\end{a1}
\em Proof of Lemma \ref{LL5.5}. \rm To simplify notation let 
$\Lambda^A_{\EI}=\exp\lt(-\frac{A}{\EI^n}\rt)$. Note by (\ref{eqeq100})
\begin{equation}
\label{band0400}
B_{\frac{r}{2}}(x)\subset \bigcup_{x\in \Xi\cap B_{\frac{r}{2}}(x)} 
B_{\frac{\CI_0}{8}r \Lambda^A_{\EI}}(x).
\end{equation}
So by Theorem 2.7 \cite{mat} we can extract some finite collection 
$$
\lt\{B_{\frac{\CI_0}{8}r \Lambda^A_{\EI}}(x_1),
B_{\frac{\CI_0}{8}r \Lambda^A_{\EI}}(x_2),\dots   B_{\frac{\CI_0}{8}r \Lambda^A_{\EI}}(x_P)\rt\}
$$ 
where 
\begin{equation}
\label{eband01}
B_{\frac{r}{2}}(x)\subset \bigcup_{i=1}^P B_{\frac{\CI_0}{8}r \Lambda^A_{\EI}}(x_i) 
\end{equation}
and 
\begin{equation}
\label{eqeq120}
\sum_{i=1}^P \cha_{B_{\frac{\CI_0}{8}r \Lambda^A_{\EI}}(x_i)}\leq c.
\end{equation}
Now for any $i,j\in \lt\{1,2,\dots P\rt\}$ if we have $B_{\frac{\CI_0}{8}r \Lambda^A_{\EI}}(x_i)\cap B_{\frac{\CI_0}{8}r \Lambda^A_{\EI}}(x_j)\not =\emptyset$ then 
\begin{equation}
\label{eqeq121}
B_{\frac{\CI_0}{8}r \Lambda^A_{\EI}}(x_i)\subset B_{\frac{\CI_0}{4}r \Lambda^A_{\EI}}(x_j). 
\end{equation}
We assume we order the balls such that $B_{\frac{\CI_0}{8}r \Lambda^A_{\EI}}(x_{i+1})\cap B_{\frac{\CI_0}{8}r \Lambda^A_{\EI}}(x_{i})\not=\emptyset$ for each 
$i=1,2,\dots P-1$. Since $\Xint{-}_{B_{\frac{r}{2}}(x_i)} \lt|S(\na u)-S(\na v)\rt|^p dz\leq 2^n \ep$, 
$\Xint{-}_{B_{\frac{r}{2}}(x_i)} \lt|\mathrm{sgn}(\det(\na v))-1\rt| dz\leq 2^n \ep$, 
$\Xint{-}_{B_{\frac{r}{2}}(x_i)} \lt|\na u\rt|^q dz\leq 2^n A$. By applying Lemma \ref{LL5} on each ball $B_{\frac{r}{2}}(x_i)$ we have that for each $i\in \lt\{1,2,\dots P\rt\}$ there 
exists $R_i\in SO(n)$ such that 
\begin{equation}
\label{eqeq125}
\Xint{-}_{B_{\frac{\CI_0}{2} \Lambda^A_{\EI} r }(x_i)} \lt|\na v-R_i\na u\rt| dz\leq 
\CI_1 A \exp\lt(\frac{2^n n A}{\EI^n}\rt)\lt(\mathrm{In}(2+\frac{\ep^{-\frac{1}{2}}}{2^n})\rt)^{-\frac{1}{32n}}.
\end{equation}
Now by (\ref{eqeq121}) $\lt|B_{\frac{\CI_0}{4}r \Lambda^A_{\EI}}(x_{i+1})\cap B_{\frac{\CI_0}{4}r \Lambda^A_{\EI}}(x_{i})\rt|\geq \frac{\CI_0^n}{8^n}  r^n (\Lambda^A_{\EI})^n$ and so by (\ref{eqeq100}) there must exists 
$$
\omega_0\in B_{\frac{\CI_0}{4}r \Lambda^A_{\EI}}(x_{i+1})\cap B_{\frac{\CI_0}{4}r \Lambda^A_{\EI}}(x_{i})\cap \Xi.
$$
So as 
\begin{equation}
\label{eqeq127}
B_{\frac{\CI_0}{4}r \Lambda^A_{\EI}}(\omega_0)\subset B_{\frac{\CI_0}{2}r \Lambda^A_{\EI}}(x_i)\cap 
B_{\frac{\CI_0}{2}r \Lambda^A_{\EI}}(x_{i+1})
\end{equation}
by definition of $\Xi$ we have that 
\begin{equation}
\label{ft1}
B_{\frac{\EI \CI_0}{8}r \Lambda^A_{\EI}}(u(\omega_0))\subset 
u\lt(B_{\frac{\CI_0}{32}r \Lambda^A_{\EI}}(\omega_0)\rt).
\end{equation}
Let $C_{iso}$ denote the constant of the isoperimetric inequality in $\R^n$. 
We claim (\ref{ft1}) implies 
\begin{equation}
\label{band0122}
\int_{B_{\frac{\CI_0}{4}r \Lambda^A_{\EI}(\omega_0)}} 
\lt|\na u\rt|^{n-1} dz\geq \frac{\EI^{n-1} \CI_0^n}{n^3 C_{iso}128^{n-1}}r^n (\Lambda_{\EI}^A)^n.
\end{equation}

Suppose this is not true. Define 
$\psi:u\lt(B_{\frac{\CI_0}{4} r \Lambda^A_{\EI}}(\omega_0)\backslash B_{\frac{\CI_0}{8} r \Lambda^A_{\EI}}(\omega_0)\rt)\rightarrow \R$ by $\psi(z)=\lt|u^{-1}\rt|$. Since by Theorem 4.1 \cite{henclkoma}
we know $u^{-1}\in W^{1,n}$ we know that $\psi\in W^{1,n}$. So either by considering difference quotients or by applying the 
Chain rule of \cite{ambdal} we have
\begin{equation}
\label{zband035}
\lt|\na \psi(z)\rt|\leq \| \na u(u^{-1}(z))^{-1}\|=
\|\mathrm{ADJ}(\na u(u^{-1}(z)))\|\det(\na u^{-1}(z)).
\end{equation}
 So 
by the Co-area formula we have 
\begin{eqnarray}
\int_{\frac{\CI_0}{8}r \Lambda^A_{\EI}}^{\frac{\CI_0}{4}r \Lambda^A_{\EI}} H^{n-1}(\psi^{-1}(t)) dt
&=&\int_{ u\lt(B_{\frac{\CI_0}{4} r \Lambda^A_{\EI}}(\omega_0)\backslash B_{\frac{\CI_0}{8} r \Lambda^A_{\EI}}(\omega_0)\rt)} \lt|\na \psi(z) \rt| dz\nn\\
&\overset{(\ref{zband035})}{\leq}& n\int_{ u\lt(B_{\frac{\CI_0}{4} r \Lambda^A_{\EI}}(\omega_0)\backslash B_{\frac{\CI_0}{8} r \Lambda^A_{\EI}}(\omega_0)\rt)} 
\lt|\mathrm{ADJ}(\na u(u^{-1}(z)))\rt| \det(\na u^{-1}(z)) dz\nn\\
&=&n\int_{B_{\frac{\CI_0}{4} r \Lambda^A_{\EI}}(\omega_0)\backslash B_{\frac{\CI_0}{8} r \Lambda^A_{\EI}}(\omega_0)} 
\lt|\mathrm{ADJ}(\na u(y))\rt| dy\nn\\
&\leq& n^3\int_{B_{\frac{\CI_0}{4} r \Lambda^A_{\EI}}(\omega_0)} \lt|\na u(y)\rt|^{n-1} dy.\nn
\end{eqnarray}
So since we are assuming (\ref{band0122}) is false 
there must exists $t\in \lt(\frac{\CI_0}{8}r \Lambda^A_{\EI} , \frac{\CI_0}{4}r \Lambda^A_{\EI}\rt)$ 
such that 
\begin{equation}
\label{zband01}
H^{n-1}\lt(\psi^{-1}(t)\rt)\leq \frac{8(\EI \CI_0)^{n-1}}{128^{n-1} C_{iso}} (r \Lambda^A_{\EI})^{n-1}.
\end{equation} 
However by construction $\psi^{-1}(t)=u(\partial B_t(\omega_0))$. Now note that by the isoperimetric 
inequality we have 
\begin{eqnarray}
\label{band0601}
\lt|u(B_t(\omega_0))\rt|^{\frac{n-1}{n}}&\leq& C_{iso} H^{n-1}(\partial u(B_t(\omega_0)))\nn\\
&\overset{(\ref{zband01})}{\leq}& \frac{8(\EI \CI_0)^{n-1}}{128^{n-1}}(r \Lambda^A_{\EI})^{n-1}.
\end{eqnarray}
Hence 
\begin{eqnarray}
\EI^n \frac{\CI_0^n}{8^n}r^n (\Lambda_{\EI}^A)^n&\overset{(\ref{ft1})}{\leq}&
\lt|u\lt(B_{\frac{\CI_0}{8} r \Lambda^A_{\EI}}(\omega_0)\rt)\rt|\nn\\
&\overset{(\ref{band0601})}{\leq}& 8^{\frac{n}{n-1}}\frac{(\EI \CI_0)^n}{128^n} (r \Lambda^A_{\EI})^n
\end{eqnarray}
which is a contradiction. Thus (\ref{band0122}) is established. 

So by (\ref{eqeq125}) and (\ref{eqeq127}) we have 
\begin{eqnarray}
\label{band0121.5}
&~&\CI_1 A \exp\lt(\frac{2^n n A}{\EI^n}\rt)\lt(\mathrm{In}(2+\frac{\ep^{-\frac{1}{2}}}{2^n})\rt)^{-\frac{1}{32n}} (\Lambda^A_{\EI})^n r^n \CI_0^n\nn\\
&~&\qd\qd\geq\int_{B_{\frac{\CI_0}{4}r \Lambda^A_{\EI}}(\omega_0)} \lt|\lt(R_i-R_{i+1}\rt)\na u\rt| dz.
\end{eqnarray}
Let 
$$
\OI=\lt\{z\in B_{\frac{\CI_0}{4} r \Lambda^A_{\EI}}(\omega_0):K(z)<\lt(\mathrm{In}(2+\frac{\ep^{-\frac{1}{2}}}{2^n})\rt)^{\frac{1}{64n}}\rt\}.
$$
So 
\begin{equation}
\label{ffeqq1}
\lt|B_{\frac{\CI_0}{4} r \Lambda^A_{\EI}}(\omega_0)\backslash \OI\rt|\leq c\lt(\mathrm{In}(2+\frac{\ep^{-\frac{1}{2}}}{2^n})\rt)^{-\frac{1}{64n}}. 
\end{equation}
Hence 
\begin{eqnarray}
\label{band0122.7}
\int_{\OI} \lt|\na u\rt|^{n-1} dz&\geq& \int_{B_{\frac{\CI_0}{4} r \Lambda^A_{\EI}}(\omega_0)} \lt|\na u\rt|^{n-1} dz
-\int_{B_{\frac{\CI_0}{4} r \Lambda^A_{\EI}}(\omega_0)\backslash \OI} \lt|\na u\rt|^{n-1} dz\nn\\
&\overset{(\ref{band0122})}{\geq}& \frac{\EI^{n-1} \CI_0^n}{n^3 C_{iso}128^{n-1}}r^n (\Lambda_{\EI}^A)^n
-\lt(\int_{B_{\frac{\CI_0}{4} r \Lambda^A_{\EI}}(\omega_0)} \lt|\na u\rt|^q dz\rt)^{\frac{n-1}{q}}\lt|B_{\frac{\CI_0}{4} r \Lambda^A_{\EI}}(\omega_0)\backslash \OI\rt|^{\frac{q-1}{q}}\nn\\
&\overset{(\ref{ffeqq1})}{\geq}&\frac{\EI^{n-1} \CI_0^n}{n^3 C_{iso}256^{n-1}}r^n (\Lambda_{\EI}^A)^n
\end{eqnarray}
And by (\ref{lband02}) 
$$
\lt|(R_i-R_{i+1})\na u(z)\rt|\geq n^{-2}\lt(\mathrm{In}(2+\frac{\ep^{-\frac{1}{2}}}{2^n})\rt)^{-\frac{1}{64n}}\lt|(R_i-R_{i+1})\rt|\lt|\na u(z)\rt|\text{ for any }z\in \OI.
$$
Putting this together with (\ref{band0121.5}) we have 
\begin{eqnarray}
\label{band0121}
&~&n^2 \CI_1 A \exp\lt(\frac{2^n n A}{\EI^n}\rt)\lt(\mathrm{In}(2+\frac{\ep^{-\frac{1}{2}}}{2^n})\rt)^{-\frac{1}{64n}} (\Lambda^A_{\EI})^n r^n \CI_0^n\nn\\
&~&\qd\qd\geq\int_{\OI} \lt|R_i-R_{i+1}\rt|\lt|\na u\rt| dz.
\end{eqnarray}
%
%
%
%
Now by (\ref{zband02}) $\int_{B_{\frac{\CI_0}{4}r \Lambda^A_{\EI}}(\omega_0)} \lt|R_i-R_{i+1}\rt|^q\lt|\na u\rt|^q dz\leq  n^{2q} 
A r^n$. Let 
$\theta=\frac{q(n-2)}{(n-1)(q-1)}$, so $\frac{1}{n-1}=\frac{\theta}{q}+1-\theta$. Note $1-\theta=\frac{q-n+1}{(n-1)(q-1)}$. 
Now letting $\tau(q)=\frac{q-n+1}{(n-1)(q-1)}$ we have $\tau'(q)=\frac{n-2}{(n-1)(q-1)^2}\geq 0$ for 
all $q\geq n$. So $(1-\theta)\geq \frac{1}{(n-1)^2}$ for all $q\geq n$. 
By the $L^p$ interpolation inequality (see 
Appendix B2 \cite{evanspde}) since $r^{\theta\frac{n}{q}}r^{n(1-\theta)}=r^{\frac{n}{n-1}}$ we know 
%
%
%
%
\begin{eqnarray}
\label{band0101}
&~&\|(R_i-R_{i+1}) \lt|\na u\rt|\|_{L^{n-1}(\OI)}\nn\\
&~&\qd\qd\leq  
\|(R_i-R_{i+1}) \lt|\na u\rt|\|_{L^{q}(B_{\frac{\CI_0}{4}r \Lambda^A_{\EI}}(\omega_0))}^{\theta}
\|(R_i-R_{i+1}) \lt|\na u\rt|\|_{L^{1}(\OI)}^{1-\theta}\nn\\
&~&\qd\qd\leq c A r^{\theta\frac{n}{q}} \|\lt(R_i-R_{i+1}\rt) \lt|\na u\rt|\|_{L^{1}(\OI)}^{1-\theta}\nn\\
&~&\qd\qd\overset{(\ref{band0121})}{\leq}   c A^2 r^{\frac{n}{n-1}} 
\exp\lt(\frac{2^n n^2 A}{\EI^n}\rt)\lt(\mathrm{In}(2+\frac{\ep^{-\frac{1}{2}}}{2^n})\rt)^{-\frac{1}{64(n-1)^2 n}} . 
\end{eqnarray}
So 
\begin{eqnarray}
\lt|R_i-R_{i+1}\rt|\frac{\EI^{n-1} \CI_0^n r^n (\Lambda_{\EI}^A)^n}{256^{n-1} n^3 C_{iso}}&\overset{(\ref{band0122.7})}{\leq}&\lt|R_i-R_{i+1}\rt|\int_{\OI} \lt|\na u\rt|^{n-1} dz\nn\\
&\overset{(\ref{band0101})}{\leq}&  c r^{n} A^{2n} \exp\lt(\frac{2^n n^3 A}{\EI^n}\rt)
\lt(\mathrm{In}(2+\frac{\ep^{-\frac{1}{2}}}{2^n})\rt)^{-\frac{1}{64 n^2}}.\nn
\end{eqnarray}
Hence 
$$
\lt|R_i-R_{i+1}\rt|\leq  c \EI^{-n} A^{2n}  
\exp\lt(\frac{2^{n+1} n^3 A}{\EI^n}\rt) \lt(\mathrm{In}(2+\frac{\ep^{-\frac{1}{2}}}{2^n})\rt)^{-\frac{1}{64n^2}}.  \nn\\
$$
Now since from (\ref{eqeq120}) 
\begin{equation}
\label{band0641}
P\leq c \lt(\Lambda^A_{\EI}\rt)^{-n}.
\end{equation}
 we know 
\begin{equation}
\label{band0350}
\lt|R_1-R_i\rt|\leq  c\EI^{-n} A^{2n}
\exp\lt(\frac{2^{n+2} n^3 A}{\EI^n}\rt)\lt(\mathrm{In}(2+\frac{\ep^{-\frac{1}{2}}}{2^n})\rt)^{-\frac{1}{64n^2}}  .
\end{equation}

So 
\begin{eqnarray}
\int_{B_{\frac{r}{2}}(x)} \lt|\na v-R_1 \na u\rt| dz&\overset{(\ref{eband01})}{\leq}&
\sum_{i=1}^P \int_{B_{\frac{\CI_0}{8} r \Lambda_{\EI}^A}(x_i)} \lt|\na v-R_1 \na u\rt| dz\nn\\
&\leq&\sum_{i=1}^P \int_{B_{\frac{\CI_0}{8} r \Lambda_{\EI}^A}(x_i)} \lt|\na v-R_i \na u\rt| dz
+\lt|(R_1-R_i)\na u\rt| dz\nn\\
&\overset{(\ref{eqeq125}),(\ref{band0641}),(\ref{band0350})}{\leq}& 
c A\exp\lt(\frac{2^n n A}{\EI^n}\rt)\lt(\mathrm{In}(2+\frac{\ep^{-\frac{1}{2}}}{2^n})\rt)^{-\frac{1}{32n}}  r^n\nn\\
&~&+c\EI^{-n} A^{2 n}\exp\lt(\frac{2^{2n+2} n^3 A}{\EI^n}\rt) \lt(\mathrm{In}(2+\frac{\ep^{-\frac{1}{2}}}{2^{n+2}})\rt)^{-\frac{1}{64n^2}}     \int_{B_r(x)} \lt|\na u\rt| dz  \nn\\
&\overset{(\ref{zband02})}{\leq}& c\EI^{-n}A^{3n} \exp\lt(\frac{2^{n+2} n^3 A}{\EI^n}\rt)  \lt(\mathrm{In}(2+\frac{\ep^{-\frac{1}{2}}}{2^n})\rt)^{-\frac{1}{64n^2}} r^n.
\end{eqnarray}

%
%
%

\begin{a1} 
\label{LL9}
Suppose we have measurable $K:\Omega\rightarrow \R_{+}$ and $u_k\in W^{1,n}(\Omega:\R^n)$ is an equibounded sequence with 
\begin{equation}
  \|\na u_k\|^n \leq K \det(\na u_k)\text{ for all }k
\end{equation}
and $(u_k)$ converges weakly in $W^{1,n}$ to $u$. For a.e.\ $x\in \Omega$ there exists 
$r_x>0$ and $N_x\in \N$ such that $u \lfloor B_{r_x}(x)$ and $u_k \lfloor B_{r_x}(x)$ are injective for all $k\geq N_x$ 
\begin{equation}
\label{opxx20} 
B_{r_x \det\lt(\na u(x)\rt)^{\frac{1}{n}}/2 K(x)^{\frac{n-1}{n}}}(u(x))\subset u\lt(B_{r_x}(x)\rt)\cap u_k(B_{r_x}(x)),
\end{equation}
\end{a1}
\em Proof of Lemma \ref{LL9}. \rm  We know that for a.e.\ $x\in \Omega$, $\det(\na u(x))>0$ 
and by Theorem 1.4. \cite{geriw} we have that $u$ is a quasiregular and satisfies $\|\na u(x)\|^n\leq 
K(x) \det(\na u(x))$. For any matrix $A\in \R^{n\times n}$ let 
$\Lambda(A)=\inf_{v\in S^{n-1}} \lt|A v\rt|$. By (\ref{fqq1})
\begin{equation}
\label{opxz1}
\Lambda\lt(\na u(x)\rt)\geq \frac{\det\lt(\na u(x)\rt)^{\frac{1}{n}}}{K^{\frac{n-1}{n}}(x)}.
\end{equation}
Pick $x$ for which $\na u(x)$ exists and 
$\det(\na u(x))>0$. Let $\delta=\frac{(\det(\na u(x)))^{\frac{1}{n}}}{100 K(x)^{\frac{n-1}{n}}}$. Let $L(z)=\na u(x)(z-x)+u(x)$. 
Let $\SI(x)=\min\lt\{\|\na u(x)\|,\|\na u(x)\|^{-1}\rt\}$. We can find $\tau_x>0$ such that 
\begin{equation}
\label{nmm1}
\lt|u(z)-L(z)\rt|\leq \frac{\delta^3}{2}\SI(x)\lt|z-x\rt|\text{ for }z\in B_{\tau_x}(x).
\end{equation}
As we have seen before in Lemma \ref{LL5}, \cite{manfredi} we know 
that for any compact subset $\wt{\Omega}\subset\subset \Omega$, letting $d(\wt{\Omega},\partial \Omega)=\sigma$ 
we have 
\begin{eqnarray}
\label{dband01}
\mathrm{osc}_{B_h(y)} u_k &\leq& c\lt(\log\lt(\frac{\sigma}{h}\rt)\rt)^{-\frac{1}{n}}
\lt(\int_{B_{\sigma}(y)} \lt|\na u_k\rt|^n\rt)^{\frac{1}{n}}\nn\\
&\leq& c\lt(\log\lt(\frac{\sigma}{h}\rt)\rt)^{-\frac{1}{n}}\text{ for any }k.
\end{eqnarray}
Hence the sequence is equi-continuous and 
\begin{equation}
\label{vbbeq1}
u_k\overset{L^{\infty}(\wt{\Omega})}{\rightarrow} u\text{ for any }\wt{\Omega}\subset\subset \Omega.
\end{equation}
So we can find $N_x\in \N$ such that for every $k\geq N_x$, 
\begin{equation}
\label{nmm2}
\lt|u_k(z)-L(z)\rt|\leq \delta^3 \SI(x)\lt|z-x\rt|  \text{ for }z\in B_{\tau_x}(x).
\end{equation}
Note that since 
$\Lambda(\na u(x))\overset{(\ref{opxz1})}{\geq} 100\delta$. So $B_{100\delta \tau_x}(L(x))\subset L(B_{\tau_x}(x))$
Let $\Psi_t(z)=t L(z)+(1-t)u_k(z)$ so by (\ref{nmm2}) 
\begin{equation}
\label{nmm4}
\Psi_t(\partial B_{\tau_x}(x))\subset N_{\delta^3 \tau_x}(L(\partial B_{\tau_x}(x)))\text{ for every }t\in \lt[0,1\rt].
\end{equation}
And in particular 
\begin{equation}
\label{nmm5}
\Psi_t(\partial B_{\tau_x}(x))\cap B_{99\delta \tau_x}(u(x))=\emptyset\text{ for every }t\in \lt[0,1\rt].
\end{equation}
Thus $\deg(u_k,B_{\tau_x}(x),z)=\deg(L,B_{\tau_x}(x),z)=1\text{ for every }t\in \lt[0,1\rt], z\in B_{99\delta \tau_x}(u(x))$. 

Now $\|\na L\|=\|\na u(x)\|$ so $L(B_{\frac{99\delta \SI(x)\tau_x}{2}}(x))\subset 
L( B_{\frac{99\delta\tau_x \|\na u(x)\|^{-1}}{2}}(x))\subset B_{\frac{99 \delta \tau_x}{2}}(u(x))$. So 
by (\ref{nmm2}) we know that for $k\geq N_x$
$$
u_k\lt(\partial B_{\frac{99\delta \SI(x)\tau_x}{2}}(x)\rt)\subset 
N_{\delta^3 \SI(x) \tau_x}\lt(L\lt(\partial B_{\frac{99 \delta\SI(x)\tau_x}{2}}(x)\rt) \rt)
\subset B_{99\delta \tau_x}(u(x)),
$$
thus $u_k\lfloor  B_{\frac{99\delta \SI(x)\tau_x}{2}}(x)$ is injective. By the same argument 
$u\lfloor  B_{\frac{99 \delta\SI(x)\tau_x}{2}}(x)$ is injective so defining  
$r_x=\frac{99 \delta\SI(x)\tau_x}{2}$ establishes the first part of the lemma. 

Now by definition of $\Lambda$ we know $B_{\Lambda(\na u(x)) r_x}(u(x))\subset L(B_{r_x}(0))$ 
and by (\ref{opxz1}) and definition of $\delta>0$ we have $\Lambda(\na u(x))\geq 100\delta$ and 
so $\delta^3\SI(x)\leq \lt(\frac{\Lambda(\na u(x))}{100}\rt)^3$, thus by 
(\ref{nmm1}), (\ref{nmm2}) $B_{\frac{\Lambda(\na u(x)) r_x}{2}}(u(x))\subset u(B_{r_x}(x))\cap u_k(B_{r_x}(x))$ hence by 
(\ref{opxz1}), (\ref{opxx20}) follows.

%
%
%
%
{\begin{a1}
\label{LL6}
Let $p\in \lt[1,n\rt]$, $q=\frac{p(n-1)}{p-1}$. Suppose $(v_k)$ is an equibounded sequence in $W^{1,p}(\Omega:\R^n)$ and $(u_k)$ an equibounded sequence in $W^{1,q}(\Omega:\R^n)$. Let $K:\Omega\rightarrow \R_{+}$ be a measurable function and assume sequence $(u_k)$
satisfies $\det(\na u_k(x))>0$ and $\|\na u_k(x)\|^n \leq K(x)\det(\na u_k(x))$ for a.e.\ $x\in \Omega$, for 
any $k\in \mathbb{N}$. Assume also 
that $\mathrm{sgn}(\det(\na v_k))\overset{L^1}{\rightarrow} 1$, 
\begin{equation}
\label{gga1}
\int_{\Omega} \lt|S(\na u_k)-S(\na v_k)\rt|^p dz\rightarrow 0\text{ as }k\rightarrow \infty
\end{equation}
and $u_k\overset{W^{1,n}}{\rightharpoonup} u$, $v_k\overset{W^{1,1}}{\rightharpoonup} v$. Then 
for a.e.\ $x\in \Omega$, $\det(\na u(x))>0$ and there exists $R_x\in SO(n)$ such that $R_x \na v(x)=\na u(x)$. Consequently 
$S(\na u(x))=S(\na v(x))$ and $\det(\na v(x))>0$ for a.e.\ $x\in \Omega$. 
\end{a1}
%
%

\em Proof of Lemma \ref{LL6}. \rm


\em Step 1. \rm Let $\omega>0$. For a.e.\ $x\in \Omega$ there exists $w_x>0$ such 
that for any $\tau\in (0,w_x)$ we can find $N_{\tau}\in \mathbb{N}$ with the property that if 
$k\geq N_{\tau}$ then for $R_k\in SO(n)$ we have  
\begin{equation}
\label{nmm5}
\Xint{-}_{B_{\tau}(x)} \lt|\na u_k-R_k \na v_k\rt| dz\leq \omega. 
\end{equation}

\em Proof of Step 1. \rm By Lemma \ref{LL9} for a.e. $x\in \Omega$ there exists $r_x>0$, 
$N_x\in \N$ such that for every $k>N_x$, $u\lfloor B_{r_x}(x)$ and $u_k\lfloor B_{r_x}(x)$ 
are injective and 
\begin{equation}
\label{mfm1}
B_{r_x\det(\na u(x))^{\frac{1}{n}}/8 K(x)^{\frac{n-1}{n}}}(u(x))\subset u(B_{\frac{r_x}{4}}(x))\cap u_k(B_{\frac{r_x}{4}}(x)).
\end{equation}
So by (\ref{vbbeq1}) we can assume $N_x$ was choosen large enough so that 
\begin{equation}
\label{mfm1.5}
B_{r_x\det(\na u(x))^{\frac{1}{n}}/16 K(x)^{\frac{n-1}{n}}}(u_k(x))\subset u(B_{\frac{r_x}{4}}(x))\cap u_k(B_{\frac{r_x}{4}}(x)).
\end{equation}
Let $\EI_x=\min\lt\{1,\frac{\det(\na u(x))^{\frac{1}{n}}}{16 K(x)^{\frac{n-1}{n}}}\rt\}$. Now $u_k$ is equibounded in $W^{1,n}$ so 
let $C_1$ be such that 
$$
\sup_{k}\int_{\Omega} \lt|\na u_k\rt|^n dz\leq C_1.
$$
Let 
$$
\tau_m=\max\lt\{\Xint{-}_{B_{r_x}(x)}\lt|S(\na u_m)-S(\na v_m)\rt|^p dz, 
\Xint{-}_{B_{r_x}(x)}\lt|\mathrm{sgn}(\det(\na v_m))-1\rt| dz\rt\}
$$
so $\tau_m\rightarrow 0$ as $m\rightarrow \infty$. So applying Lemma \ref{LL5} on $B_{r_x}(x)$ we have that 
for 
$A_x:=\frac{C_1}{r_x^n}$ we have for some $R_m\in SO(n)$
$$
\Xint{-}_{B_{\CI_0 r_x \exp\lt(-\frac{A_x}{\EI_x}\rt)}(x)} \lt|\na u_m-R_m \na v_m\rt| dz\leq \CI_1 A_x
\exp\lt(\frac{n A_x}{\EI_x}\rt) \mathrm{In}\lt(2+\tau_m^{-\frac{1}{4}}\rt)^{-\frac{1}{32n}}.
$$
So defining $w_x=\CI_0 r_x \exp\lt(-\frac{A_x}{\EI_x}\rt)$. For any $\tau\in (0,w_x)$ we can find 
$N_{\tau}$ such that for $k\geq N_{\tau}$ for some $R_k\in SO(n)$ (\ref{nmm5}) holds true. \nl

%
%

\em Step 2. \rm  For $\sigma>0$. For a.e.\ $x\in \Omega$, $r>0$ define  
\begin{equation}
\label{ghf150}
\BI_r^{\sigma,x}:=\lt\{z\in B_r(x):\lt|u(z)-u(x)-\na u(x)(z-x)\rt|<\sigma\lt|z-x\rt|\rt\}
\end{equation}
and
\begin{equation}
\label{ghf150.5}
\DI_r^{\sigma,x}:=\lt\{z\in B_r(x):\lt|v(z)-v(x)-\na v(x)(z-x)\rt|<\sigma\lt|z-x\rt|\rt\}.
\end{equation}
Now for a.e.\ $x\in \Omega$ there exists $\mu\in (0,w_x)$ such that for any $\phi\in S^{n-1}$ we can find $y_1\in B_{\sigma \mu}(x)\cap 
\BI_{\mu}^{\sigma,x}\cap \DI_{\mu}^{\sigma,x}$ and $y_2\in A(x,\frac{\mu}{2},\mu)\cap 
\BI_{\mu}^{\sigma,x}\cap \DI_{\mu}^{\sigma,x}$ such that 
\begin{equation}
\label{fgh300}
 \lt|\lt(\frac{y_2-y_1}{\lt|y_2-y_1\rt|}\rt)-\phi\rt|\leq \sigma^{\frac{1}{n-1}}
\end{equation}
and for affine function $L_{R_x}$ with $\na L_{R_x}=R_x\in SO(n)$
\begin{equation}
\label{fgh302}
\lt|u(y_i)-L_{R_x}(v(y_i))\rt|\leq \sigma \mu\text{ for }i=1,2.
\end{equation}

In addition
\begin{equation}
\label{ghf151}
\lt|B_r(z)\backslash \BI_r^{\sigma,z}\rt|\leq \sigma^{4n} r^n\text{ and }
\lt|B_r(z)\backslash \DI_r^{\sigma,z}\rt|\leq \sigma^{4n} r^n\text{ for all }r\in \lt(0,\mu\rt].
\end{equation}

\em Proof of Step 2. \rm By Theorem 1.4 \cite{geriw} 
we know $u$ is a mapping of integrable and $\frac{\|\na u(z)\|^n}{\det(\na u(z))} \leq K(z)$ 
for a.e.\ $z\in \Omega$. So in particular $\det(\na u(z))>0$ for a.e.\ $z\in \Omega$. By Theorem 1, Section 6.1.1. \cite{evans2} for a.e.\ $x\in \Omega$ we can find 
$\mu\in (0,w_x)$ such that (\ref{ghf151}) holds true and  
\begin{equation}
\label{sula1}
\Xint{-}_{B_r(x)} \lt|\na u(x)-\na u(z)\rt| dz\leq \sigma^{4n}\text{ and }\Xint{-}_{B_r(x)} \lt|\na v(x)-\na v(z)\rt| dz\leq \sigma^{4n}\text{ for all }r\in \lt(0,\mu\rt]
\end{equation}
Fix an $x$ for which this is true and for which Step 1 holds. By Step 1 we can find $N_x\in \N$ such that 
\begin{equation}
\label{asz1.5}
\int_{B_{\mu}(x)} \lt|\na u_k- R_k\na v_k\rt| dz\leq \sigma^{4n} \frac{\mu^n}{2}\text{ for all }k\geq N_x.
\end{equation}
Passing to a subsequence (not relabeled) we have $R_k\rightarrow R_x$ as $k\rightarrow \infty$. As 
$u_k,v_k$ are equibounded in $L^1$ this implies there exists $M_x\geq N_x$ such that 
\begin{equation}
\label{asz1}
\int_{B_{\mu}(x)} \lt|\na u_k- R_x\na v_k\rt| dz\leq \sigma^{4n} \mu^n\text{ for all } k\geq M_x.
\end{equation}

Now pick $k\geq M_x$ large enough so that 
\begin{equation}
\label{ghf153}
\|u_k-u\|_{L^{1}(\Omega)}\leq \sigma^{4n} \mu^{n+1},\; \|v_k-v\|_{L^{1}(\Omega)}\leq \sigma^{4n} \mu^{n+1}.
\end{equation}
By Poincare inequality from (\ref{asz1}) there exists affine function $L_{R_x}$ with $\na L_{R_x}=R_x$ with 
$$
\int_{B_{\mu}(x)} \lt|u_k-L_{R_x}\circ v_k\rt| dz\leq c\sigma^{4n} \mu^{n+1}. 
$$
Hence by (\ref{ghf153}) we have 
\begin{equation}
\label{band01}
\int_{B_{\mu}(x)} \lt|u-L_{R_x}\circ v\rt| dz\leq c\sigma^{4n} \mu^{n+1}. 
\end{equation}
So let 
\begin{equation}
\label{gfh200}
\GI:=\lt\{z\in B_{\mu}(x):\lt|u(z)-L_{R_x}(v(z))\rt|\leq \sigma \mu \rt\}
\end{equation}
thus 
\begin{equation}
\label{gfh201}
\lt|B_{\mu}(x)\backslash \GI\rt|\leq c\sigma^{4n-1} \mu^n.
\end{equation}
Let $\HI:=\GI\cap \lt(\BI_{\mu}^{\sigma,x} \cup\DI_{\mu}^{\sigma,x}\rt)$. Since 
by (\ref{ghf151}) and (\ref{gfh201})
\begin{equation}
\label{band090}
\lt|B_{\mu}(x)\backslash \HI\rt|\leq c\sigma^{4n-1} \mu^n. 
\end{equation}
Now let 
\begin{equation}
\label{band091}
\EI:=\lt\{y\in B_{\sigma \mu}(x):\int \cha_{B_{\mu}(x)\backslash \HI}(z)\lt|z-y\rt|^{-n+1} dz\leq \sigma^2 \mu\rt\}
\end{equation}
Now 
\begin{eqnarray}
\label{cct1}
&~&\int_{B_{\sigma \mu}(x)} \int \cha_{B_{\mu}(x)\backslash \HI}(z)\lt|z-y\rt|^{-n+1} dz dy\nn\\
&~&\qd\qd\qd=\int \cha_{B_{\mu}(x)\backslash \HI}(z)  \lt( \int_{B_{\sigma \mu}(x)} \lt|z-y\rt|^{-n+1} dy\rt) dz\nn\\
&~&\qd\qd\qd =c\sigma\mu\int \cha_{B_{\mu}(x)\backslash \HI} (z)  dz\nn\\
&~&\qd\qd\qd \overset{(\ref{band090})}{\leq} c\sigma^{4n}\mu^{n+1}.
\end{eqnarray}
Hence $\lt|B_{\sigma\mu}(x)\backslash \EI\rt|\sigma^2 \mu\overset{(\ref{band091}),(\ref{cct1})}{\leq} 
c\sigma^{4n} \mu^{n+1}$ so 
\begin{equation}
\label{band06}
\lt|B_{\sigma \mu}(x)\backslash \EI\rt|\leq c\sigma^{4n-2} \mu^n.
\end{equation}

By (\ref{band090}), (\ref{band06}) it is clear 
$\lt|B_{\sigma \mu}(x)\cap (\EI\cup \HI)\rt|>0$ so pick $y_1\in \EI\cap \HI\cap B_{\sigma \mu}(x)$, define $l_{y_1}^{\theta}:=\lt\{y_1+\R_{+}\theta\rt\}$. 
So by the Co-area formula into $S^{n-1}$ (see for example \cite{JL} Lemma 14), by definition of $\EI$ (recall (\ref{band091})) 
\begin{equation}
\label{gfh230}
\int_{\theta\in S^{n-1}} \int_{l_{y_1}^{\theta}} \cha_{B_{\mu}(x)\backslash \HI} dH^{1} z dH^{n-1} \theta\leq 
\sigma^2 \mu.
\end{equation}
So let 
$$
\Psi_{y_1}:=\lt\{\theta\in S^{n-1}: \int_{l_{y_1}^{\theta}} \cha_{B_{\mu}(x)\backslash \HI} dH^1 z\leq 
\sigma\mu\rt\}
$$ 
thus by (\ref{gfh230}), $H^{n-1}(S^{n-1}\backslash \Psi_{y_1})\leq \sigma$. Thus we can find 
$\psi\in \Psi_{y_1}$ such that $\lt|\psi-\theta\rt|\leq c\sigma^{\frac{1}{n-1}}$. Thus we can find $y_2\in 
l_{y_1}^{\psi}\cap A(x,\frac{\mu}{2},\mu)\cap \HI$. Since 
$y_1,y_2\in \HI\subset \GI$ by definition (\ref{gfh200}) we know they satisfy (\ref{fgh302}). Since 
$\frac{y_2-y_1}{\lt|y_2-y_1\rt|}=\psi$ it is clear that (\ref{fgh300}) is satisfied. This completes the proof of Step 2. \nl

\em Step 3. \rm We will show that for a.e.\ $x\in \Omega$ there exists 
$R_x\in SO(n)$ such that 
\begin{equation}
\label{band08}
R_x \na v(x)=\na u(x).
\end{equation}

\em Proof of Step 3. \rm Let $x\in \Omega$ be one of the a.e.\ points $x$ 
such that the conclusion of Step 2 hold true. 
Let $\gamma>0$ and set 
\begin{equation}
\label{jk12}
\sigma=\lt(\frac{\gamma}{\lt|\na u(x)\rt|+\lt|\na v(x)\rt|+1}\rt)^{n-1}. 
\end{equation}
By Step 3 we can find and points $y_1\in B_{\sigma \mu}(x)$, 
$y_2\in A(x,\frac{\mu}{2},\mu)$ such that (\ref{fgh300}),  (\ref{fgh302}) are satisfied. So since $y_i\in 
\BI^{\sigma,x}_{\mu}$
$$
\lt|u(y_i)-u(x)-(y_i-x)\cdot \na u(x)\rt|<\sigma \lt|y_i-x\rt|<\sigma \mu\text{ for }i=1,2
$$
taking one inequality away from another 
\begin{equation}
\label{gfh140}
\lt|(u(y_1)-u(y_2))-\na u(x)(y_1-y_2)\rt|<2\sigma \mu.
\end{equation}
And in the same way 
\begin{equation}
\label{gfh141}
\lt|(v(y_1)-v(y_2))-\na v(x)(y_1-y_2)\rt|<2\sigma \mu.
\end{equation}

Applying (\ref{fgh302}) to (\ref{gfh140}) we have 
$\lt|R_x(v(y_1)-v(y_2))-\na u(x)(y_1-y_2)\rt|<4\sigma \mu$ 
and putting this together with (\ref{gfh141}) we 
$\lt| \na v(x)(y_1-y_2)-R_x^{-1}\na u(x)(y_1-y_2)\rt|\leq 6\sigma \mu$ since $\lt|y_1-y_2\rt|>\frac{\mu}{2}$ so 
\begin{equation}
\label{jkl1}
\lt| \na v(x)\frac{(y_2-y_1)}{\lt|y_2-y_1\rt|}-R_x^{-1}\na u(x)\frac{(y_2-y_1)}{\lt|y_2-y_1\rt|}\rt|\leq 12\sigma. 
\end{equation}
Thus $\lt|(\na v(x)-R_x^{-1}\na u(x))\phi\rt|\overset{(\ref{jkl1})}{\leq} 12\sigma+\lt|(\na v(x)-R_x^{-1}\na u(x))\rt| \lt|\phi-\frac{y_2-y_1}{\lt|y_2-y_1\rt|}\rt|\overset{(\ref{fgh300}),(\ref{jkl1}),(\ref{jk12})}{\leq} 14 \gamma$. 
Now as $\gamma$ is arbitrary this implies $R_x \na v(x)=\na u(x)$. This completes the proof of 
the lemma. $\Box$ \nl

%
%
%

\subsection{ Proof of Theorem \ref{TT2} completed}

Let $v\in W^{1,p}(\Omega:\R^n)$ and $u\in W^{1,q}(\Omega:\R^n)$ 
be the weak limit of $v_k,u_k$. We know by 
Lemma \ref{LL6} $S(\na u)=S(\na v)$ for a.e.\ $x\in \Omega$. So we 
can apply Theorem \ref{T1} and thus there exists $R\in SO(n)$ such that 
\begin{equation}
\label{xcxx271}
\na v(z)=R\na u(z)\text{ for a.e. }z\in \Omega. 
\end{equation}

Since $\det(\na u(z))>0$ for a.e.\ $z\in \Omega$ so  
$\int_{\Omega} \det(\na u) dz\leq \int_{\Omega} \|\na u\|^n dz\leq C$. For any 
$\gamma>0$ let 
\begin{equation}
\label{xcx10}
\DI_{\gamma}:=\lt\{z\in \Omega: \det(\na u(z))<\gamma^{\frac{1}{100}}\rt\}
\end{equation}
and let 
\begin{equation}
\label{xcx11}
\BI_{\gamma}:=\lt\{z\in \Omega: \lt|\na u(z)\rt|>\gamma^{-\frac{1}{100}}\rt\}.
\end{equation}
Note $\lt|\BI_{\gamma}\rt|\rightarrow 0$ and $\lt|\DI_{\gamma}\rt|\rightarrow 0$ 
as $\gamma\rightarrow 0$. Let $\delta\in (0,1)$. Define
\begin{equation}
\label{ghf500}
\OI_{\delta}:=\lt\{x\in \Omega:K(x)\geq \delta^{-\frac{1}{100}}\rt\},
\end{equation}
note 
\begin{equation}
\label{opps101.5}
\lt|\OI_{\delta}\rt|\leq c\delta^{\frac{1}{100}}. 
\end{equation}

Let $\ep\in (0,\delta)$ be small enough so that 
\begin{equation}
\label{cvczz1}
\int_{\BI_{\ep}} \det(\na u) dz\leq \delta^n.
\end{equation}

For a.e.\ $x\in \Omega\backslash \lt(\DI_{\delta}\cup \OI_{\delta}\rt)$ let $L_x$ be the affine map defined 
by $L_x(z)=u(x)+\na u(x)z$. Recall from Lemma \ref{LL9} we defined 
$\Lambda(A):=\inf_{v\in S^{n-1}} \lt|A v\rt|$ and from (\ref{fqq1}) we know 
\begin{equation}
\label{band063}
\frac{\det(\na u(x))}{K^{n-1}(x)}\leq \lt(\Lambda(\na u(x))\rt)^n\text{ for a.e. }x\in \Omega.
\end{equation}
Since $x \in \Omega\backslash \lt(\DI_{\delta}\cup \OI_{\delta}\rt)$, $K(x)< \delta^{-\frac{1}{100}}$ and $\det(\na u(x))>\delta^{\frac{1}{100}}$ 
so 
$\delta^{\frac{1}{100}}\overset{(\ref{band063})}{\leq} \Lambda(\na u(x))$. 

Thus  
\begin{equation}
\label{xcxv1}
B_{\delta^{\frac{1}{100}} h}(u(x))\subset L_{x}(B_h(x))\text{ for all }h>0.  
\end{equation}
Now by uniform continuity of $u$, approximate differentiability of $u$ and approximate continuity of $\det(\na u)$ 
for a.e.\ $x\in \Omega\backslash (\DI_{\delta}\cup \OI_{\delta})$ there exists $p_x>0$ such that 
\begin{equation}
\label{klk1}
H_d\lt(u(B_h(x)),L_x(B_h(x))\rt)\leq \ep\delta^{\frac{1}{100}} h\text{ for any }h\in (0,p_x)
\end{equation}
and 
\begin{equation}
\label{klk2}
\lt|\lt|u(B_h(x))\rt|-\Gamma(n)\det(\na u(x))h^n\rt|\leq \ep h^n\text{ for any }h\in (0,p_x).
\end{equation}
Note by (\ref{xcxv1}), (\ref{klk1}) we have that
\begin{equation}
\label{band02}
B_{8^{-1}\delta^{\frac{1}{100}}h}(u(x))\subset u(B_{\frac{h}{4}}(x))\text{ for any }h\in (0,p_x),\; x\in \Omega\backslash (\DI_{\delta}\cup \OI_{\delta}) 
\end{equation}
Now let  
\begin{equation}
\label{kjk3}
\HI_{\rho}:=\lt\{x\in \Omega\backslash (\DI_{\delta}\cup \OI_{\delta}):p_x< \rho\rt\}.
\end{equation}
Let us choose $\rho_0$ be such that 
\begin{equation}
\label{kjk4}
\lt|\HI_{\rho_0}\rt|< \delta.
\end{equation}
Note also that by Lebesgue density theorem  for a.e.\ $x\in \Omega\backslash (\DI_{\delta}\cup \OI_{\delta}\cup \HI_{\rho_0})$ the 
ratio 
$$
\frac{\lt|B_r(x)\cap (\DI_{\delta}\cup \OI_{\delta}\cup \HI_{\rho_0})\rt|}{r^n}
$$ 
can be arbitrarily small. For each $x$ we need to find the 
$q_x>0$ such that for all $r\in (0,q_x)$ the ratio is less than a small constant depending on $\delta$. Rather than introduce more 
notation to signify this small quantity then later take it to be less than the constant we need, we find $q_x$ that has the 
exact property we need in terms of $\delta$. So for a.e.\ $x\in \Omega\backslash (\DI_{\delta}\cup \OI_{\delta}\cup \HI_{\rho_0})$ there exists $q_x>0$ such that 
\begin{equation}
\label{band0231}
\lt|B_r(x)\cap ( \DI_{\delta}\cup \OI_{\delta}\cup \HI_{\rho_0})\rt|\leq \frac{\CI_0^n}{32^n}r^n\exp(-64^n n\delta^{-\frac{101n}{100}})\text{ for }r\in (0,q_x).
\end{equation}
For $\rho>0$ let 
\begin{equation}
\label{band0234}
\Theta_{\rho}:=
\lt\{x\in \Omega\backslash (\DI_{\delta}\cup \OI_{\delta}\cup \HI_{\rho_0}):q_x<\rho\rt\}. 
\end{equation}
We can find $\rho_1>0$ such that 
\begin{equation}
\label{band0.018}
\lt|\Theta_{\rho_1}\rt|\leq \delta.
\end{equation}

By Lemma \ref{LL9} for a.e.\ $x\in \Omega$ there exists $w_x\in (0,q_x)$ and $N_x\in N$ such that 
for $u\lfloor B_{w_x}(x)$ and $u_k\lfloor B_{w_x}(x)$ are injective for any $k\geq N_x$. We can find $\ri>0$ and $M_0\in \N$ such that 
\begin{equation}
\label{ghf400.6}
\lt|\lt\{x\in\Omega:w_x<\tau\rt\}\rt|<\ep
\text{ and }\lt|\lt\{x\in\Omega:N_x>M_0\rt\}\rt|<\ep.
\end{equation}
Let $\EI_{\ep}:=\lt\{x\in \Omega: w_x<\tau\rt\}\cup \lt\{x\in \Omega: N_x>M_0\rt\}$. Now recall the 
notation $N_{h}(\cdot)$ (see (\ref{nna1})). Define      
\begin{equation}
\label{opps100}
\Pi:=\Omega\backslash \lt(\OI_{\delta}\cup \DI_{\delta}\cup \HI_{\rho_0}\cup \BI_{\delta}   \cup 
\EI_{\ep}\cup \Theta_{\rho_1}\cup N_{\delta}(\partial \Omega)\rt). 
\end{equation}

\em Step 1. \rm Let $\eta=\frac{1}{2}\min\lt\{\tau,\rho_0\rt\}$. 
Since $\lt\{B_{\eta}(x):x\in \Pi\rt\}$ is a cover of $\Pi$, again by Theorem 2.7 \cite{mat} we can find a 
collection $\lt\{B_{\frac{\eta}{2}}(x_1),B_{\frac{\eta}{2}}(x_2),\dots B_{\frac{\eta}{2}}(x_m)\rt\}$ such that 
\begin{equation}
\label{zband041}
\Pi\subset \bigcup_{i=1}^{m} B_{\frac{\eta}{2}}(x_i) 
\end{equation}
and 
\begin{equation}
\label{sz1}
\sum_{k=1}^m \cha_{B_{2\eta}(x_i)}\leq c.
\end{equation}
Let $\gamma$ be some small positive number we decide on later. 
Let $M_1>M_0$ be such that 
$$
\int_{\Omega}  \lt|S(\na u_q)-S(\na v_q)\rt|^p+\lt|1-\mathrm{sgn}(\det(\na v_q))\rt| dz\leq \gamma\text{ for all }q>M_1. 
$$
Fix $k>M_1$, define
\begin{equation}
\label{vvv7}
B^k_0:=\lt\{i\in\lt\{1,2,\dots m\rt\}:\Xint{-}_{B_{\eta}(x_i)} \lt|\na u_k\rt|^n  dz\geq 
\delta^{-n}\rt\}, 
\end{equation}
%
%
\begin{equation}
\label{ghf401.5}
B^k_1:=\lt\{i\in\lt\{1,2,\dots m\rt\}:\Xint{-}_{B_{\eta}(x_i)} \lt|S(\na u_k)-S(\na v_k)\rt|^p+\lt|1-\mathrm{sqn}(\det(\na v_k))\rt| dz\geq \sqrt{\gamma}\rt\}, 
\end{equation}
%
%
\begin{equation}
\label{xcx2}
B_2:=\lt\{i\in\lt\{1,2,\dots m\rt\}:
\lt|B_{\eta}(x_i)\cap \DI_{\delta}\rt|\geq \Gamma(n)2^{n-1}\eta^n \rt\}, 
\end{equation}
\begin{equation}
\label{aklk2}
B_3:=\lt\{i\in \lt\{1,2,\dots m\rt\}: \int_{B_{\eta}(x_i)\cap \BI_{\ep}} \det(\na u) dz\geq \delta^{\frac{n}{4}}\int_{B_{\eta}(x_i)} \det(\na u) dz \rt\}.
\end{equation}
We will show 
\begin{equation}
\label{cvc30}
\ca{B^k_0\cup B^k_1\cup B_2\cup B_3}\leq \lt(c\sqrt{\gamma}+c\delta^{\frac{n}{2}}+c\lt|\DI_{\delta}\rt|\rt)\eta^{-n}
\end{equation}
and any $i\in\lt\{1,2,\dots m\rt\}\backslash (B^k_0\cup B^k_1\cup B_2\cup B_3)$ there 
exists $R^k_i\in SO(n)$ such that 
\begin{equation}
\label{cvc30.7}
\Xint{-}_{B_{\frac{\eta}{2}}(x_i)} \lt|\na v_k-R^k_i \na u_k\rt| dz\leq c\ep^{2n}
\end{equation}
and (recalling $R\in SO(n)$ satisfies (\ref{xcxx271})) for affine maps $l_R,l_{R^k_i}$ with $\na l_R=R$, $\na l_{R_i^k}=R_i^k$
\begin{equation}
\label{cvc31}
\eta^{-n}\int_{u\lt(B_{\eta}(x_i)\backslash \BI_{\ep}\rt)} \lt|l_R(z)-l_{R^k_i}(z)\rt|\det(\na u^{-1}(z)) dz
\leq c\eta\ep^{2n}.
\end{equation}

\em Proof of Step 1. \rm Note that since $u_k$ is an equibounded sequence in $W^{1,n}$, so 
$$
c\delta^{-n} \eta^n \ca{B^k_0}\leq \int_{B_{\eta}(x_i)} \lt|\na u_k\rt|^n dz\leq c
$$
thus $\ca{B_0^k}\leq \frac{c \delta^n}{\eta^n}$. It is also clear
$\sqrt{\gamma}\eta^n \ca{B^k_1}\leq c\gamma$ and thus $\ca{B^k_1}\leq \sqrt{\gamma}\eta^{-n}$. Also 
$$
\ca{B_2}\eta^n\leq c\lt|\DI_{\delta}\rt|
$$ 
so $\ca{B_2}\leq c \lt|\DI_{\delta}\rt|\eta^{-n}$. And
\begin{eqnarray}
c\delta^n&\overset{(\ref{cvczz1}),(\ref{sz1})}{\geq}&\sum_{i\in B_3\backslash B_2}  
\int_{B_{2\eta}(x_i)\cap \BI_{\ep}} \det(\na u) dx\nn\\
&\overset{(\ref{aklk2})}{\geq}&
c\delta^{\frac{n}{4}}\sum_{i\in B_3\backslash B_2}   \int_{B_{2\eta}(x_i)} \det(\na u) dx\nn\\
&\overset{(\ref{xcx2}),(\ref{xcx10})}{\geq}& c\delta^{\frac{n}{2}}\eta^n \ca{B_3\backslash B_2}.\nn
\end{eqnarray}
Thus $\ca{B_3\backslash B_2}\leq  c\delta^{\frac{n}{2}}\eta^{-n}$. Now by 
definition of $\Theta_{\rho_1}$ (see (\ref{band0234}), (\ref{band0231})), since we know 
from Step 1 $\eta<\rho_0$ and $x_1,x_2,\dots x_m\not\in \Theta_{\rho_0}$ 
we have that for each $i=1,2,\dots m$
\begin{equation}
\label{vvv3}
\lt|B_{\eta}(x_i)\cap (\DI_{\delta}\cup \OI_{\delta}\cup \HI_{\rho_0})\rt|\leq 
\frac{\CI_0^n}{32^n} \eta^n \exp(-64^n n \delta^{-\frac{101n}{100}}).
\end{equation} 
And by (\ref{band02}), (\ref{kjk3})  
\begin{equation}
\label{qzqza1}
B_{8^{-1}\delta^{\frac{1}{100}}h}(u(z))\subset u(B_{\frac{h}{4}}(z))\text{ for each }
z\in B_{\frac{\eta}{2}}(x_i)\backslash (\DI_{\delta}\cup \OI_{\delta}\cup \HI_{\rho_0}), h\in (0,\frac{\eta}{2}), i=1,2,\dots m.
\end{equation}
Recalling from (\ref{dband01}) we know equicontinuity of the sequence $u_k$ on a compact 
subset of $\Omega$ and hence uniform convergence of $u_k$. So let $Q\in \N$ be such that 
\begin{equation}
\label{band0505}
\|u-u_k\|_{L^{\infty}(\Omega\backslash N_{\frac{\delta}{2}}(\partial \Omega))}<16^{-3}\eta\delta^{\frac{1}{100}}\CI_0\exp(-64^n \delta^{-\frac{101 n}{100}})\text{ for all }k\geq Q.
\end{equation}
So (recalling that $B_{\frac{\eta}{2}}(x_i)\subset N_{\frac{\eta}{2}}(\Pi)\overset{(\ref{opps100})}{\subset} \Omega\backslash N_{\frac{\delta}{2}}(\partial \Omega))$
\begin{equation}
\label{band01001}
u(B_{\frac{h}{4}}(z))\backslash N_{16^{-2}\eta\delta^{\frac{1}{100}}\CI_0\exp(-64^n \delta^{-\frac{101 n}{100}})}(\partial u(B_{\frac{h}{4}}(z)))\overset{(\ref{band0505})}{\subset} u_k(B_{\frac{h}{4}}(z))\text{ for any }k\geq Q.
\end{equation}
Thus for any $z\in B_{\frac{\eta}{2}}(x_i)\backslash (\DI_{\delta}\cup \OI_{\delta}\cup \HI_{\rho_0})$, 
$h\in \lt[\frac{\CI_0}{16}\eta \exp(-64^n  \delta^{-\frac{101n}{100}}),\frac{\eta}{2}\rt]$ we have 
\begin{eqnarray}
d(B_{16^{-1}\delta^{\frac{1}{100}}h}(u(z)),\partial u(B_{\frac{h}{4}}(z)))&\overset{(\ref{qzqza1})}{\geq}& 
16^{-1}\delta^{\frac{1}{100}}h \nn\\
&\geq & 16^{-2} \delta^{\frac{1}{100}}\CI_0 \eta \exp(-64^n \delta^{-\frac{101n}{100}})
\end{eqnarray}
and thus 
\begin{eqnarray}
\label{zband011}
B_{16^{-1}\delta^{\frac{1}{100}}h}(u(z))&\subset&
 u(B_{\frac{h}{4}}(z))\backslash N_{16^{-2} \eta\delta^{\frac{1}{100}}\CI_0 \eta \exp(-64^n \delta^{-\frac{101n}{100}})}(\partial u(B_{\frac{h}{4}}(z)))\nn\\
&\overset{(\ref{band01001})}{\subset}& u_k(B_{\frac{h}{4}}(z))\text{ for }k\geq Q.
\end{eqnarray}
Hence as $\|u-u_k\|_{L^{\infty}(\Omega\backslash N_{\frac{\delta}{2}}(\partial \Omega))}\overset{(\ref{band0505})}{\leq} 16^{-2}\delta^{\frac{1}{100}}h$ we have 
\begin{eqnarray}
\label{vv1}
B_{\frac{15}{16^2}\delta^{\frac{1}{100}}h}(u_k(z))&\overset{(\ref{zband011})}{\subset}& 
u_k(B_{\frac{h}{4}}(z))\text{ for }z\in B_{\frac{\eta}{2}}(x_i)\backslash (\DI_{\delta}\cup \OI_{\delta}\cup \HI_{\rho_0}), \nn\\
&~&\qd\qd h\in \lt[\frac{\CI_0}{16}\eta \exp(-64^n \delta^{-\frac{101n}{100}}),\frac{\eta}{2}\rt]\text{ and }k\geq Q.
\end{eqnarray}

If $i\in\lt\{1,2,\dots m\rt\}\backslash (B^k_1\cup B^k_0)$ we can define 
$\Xi:=B_{\frac{\eta}{2}}(x_i)\backslash (\DI_{\delta}\cup \OI_{\delta}\cup \HI_{\rho_0})$ and 
define 
\begin{equation}
\label{cceq1}
\EI=\frac{15}{16^2}\delta^{\frac{1}{100}}\text{ and } A=\delta^{-n}2^n 
\end{equation}
and notice that 
$$
\frac{\CI_0}{8}\frac{\eta}{2}\exp\lt(-\frac{A}{\EI^n}\rt)>\frac{\CI_0}{16} \eta 
\exp(-64^n  \delta^{-\frac{101n}{100}})\text{ and }
\frac{\CI_0^n}{16^n}\frac{\eta^n}{2^n}\exp\lt(-\frac{n A}{\EI^n}\rt)>\frac{\CI_0^n}{32^n} \eta^n 
\exp(-64^n n \delta^{-\frac{101n}{100}})
$$
thus by (\ref{vv1}), (\ref{vvv3}) hypotheses for ($r$ taken to be $\frac{\eta}{2}$) (\ref{vvv2}) and (\ref{eqeq100}) is satisfied. So we can apply Lemma \ref{LL5.5} (taking $A$, $\EI$ defined by (\ref{cceq1}) and $\ep=\sqrt{\gamma}$. In 
addition in view of (\ref{vvv7}), (\ref{ghf401.5}) hypotheses (\ref{zband03}), (\ref{zband03.5}) and 
(\ref{zband02}) are satisfied and there exists $R^k_i\in SO(n)$ such that 
\begin{equation}
\label{ghf49}
\Xint{-}_{B_{\frac{\eta}{4}}(x_i)} 
\lt|\na v_k-R^k_i \na u_k\rt| dx\leq 
c\delta^{-4 n^2} \exp\lt(40^{n+1} 5 n^3 \delta^{-\frac{101n}{100}}\rt) 
\lt(\mathrm{In}\lt(2+\frac{\gamma^{-\frac{1}{8}}}{2^n}\rt)\rt)^{-\frac{1}{64n^2}}
\end{equation}
so assuming $\gamma$ was chosen small enough (\ref{cvc30.7}) is established. 

By Poincare's inequality there exists affine map $l_{R^k_i}$ with $\na l_{R^k_i}=R^k_i$, so  
\begin{equation}
\label{ghf50.4}
\Xint{-}_{B_{\frac{\eta}{4}}(x_i)} \lt|v_k-l_{R^k_i}\circ u_k\rt| dx\leq c
\eta\ep^{2n}.
\end{equation}
Now as $v_k\overset{L^{1}(\Omega)}{\rightarrow} v$ and $u_k\overset{L^{1}(\Omega)}{\rightarrow} u$ so assuming 
$k$ is large enough we have 
$$
\Xint{-}_{B_{\frac{\eta}{4}}(x_i)} \lt|v_k-v\rt| dx\leq c\eta \ep^{2n}\text{ and } 
\Xint{-}_{B_{\frac{\eta}{4}}(x_i)} \lt|u_k-u\rt| dx\leq c\eta \ep^{2n}
$$
putting this together with (\ref{ghf50.4}) we have 
\begin{equation}
\label{ghf50}
\Xint{-}_{B_{\frac{\eta}{4}}(x_i)} \lt|v-l_{R^k_i}\circ u\rt| dx\leq c\eta\ep^{2n}.
\end{equation}

Since $\na v=R \na u$ for some affine map $l_R$ with $\na l_R=R$ we have $v=l_R\circ u$ on $\Omega$. So putting 
this together with (\ref{ghf50}) we have 
\begin{equation}
\label{ghf51}
\Xint{-}_{B_{\eta}(x_i)} \lt|l_R\circ u-l_{R^k_i}\circ u\rt| dx\leq c\eta\ep^{2n}.
\end{equation}

Thus if $i\not\in (B_0^k\cup B_1^k\cup B_2\cup B_3)$
\begin{eqnarray}
\label{fgh401}
c \eta\ep^{2n}&\geq&  \Xint{-}_{B_{\frac{\eta}{4}}(x_i)\backslash \BI_{\ep}} \lt|l_R(u(z))-l_{R^k_i}(u(z))\rt|\det((\na u(u^{-1}(u(z))))^{-1})\det(\na u(z)) dz\nn\\
&=&c \eta^{-n}\int_{u\lt(B_{\frac{\eta}{4}}(x_i)\backslash \BI_{\ep}\rt)} \lt|l_R(z)-l_{R^k_i}(z)\rt|\det(\na u^{-1}(z)) dz.
\end{eqnarray}
This completes the proof of Step 1. \nl

\em Step 2. \rm We will show that for any $k>M_1$ 
\begin{equation}
\label{vxx1}
\lt|R-R^k_i\rt|\leq c\ep^{\frac{n}{4}}\text{ for any }i\in \lt\{1,2,\dots M\rt\}\backslash (B^k_0\cup B^k_1\cup B_2\cup B_3).
\end{equation}

\em Proof of Step 2. \rm  Now since $i\not \in  B_3$, (see (\ref{aklk2}) for the definition) and we chose $x_i\not \in \BI_{\delta}$ 
(recall (\ref{xcx11}) for the definition)
\begin{eqnarray}
\label{xcxxz101}
\lt|u\lt(\BI_{\ep}\cap B_{\eta}(x_i)\rt)\rt|&\overset{(\ref{aklk2})}{<}& \delta^{\frac{n}{4}}
\lt|u\lt(B_{\eta}(x_i)\rt)\rt| \nn\\
&\overset{(\ref{klk2}),(\ref{xcx11})}{\leq}&c\delta^{\frac{n}{4}}\delta^{-\frac{n}{100}}\eta^n  \nn\\
&\leq& c\delta^{\frac{24n}{100}}\eta^n .
\end{eqnarray}
So as $x_i\in \Pi$ thus $x_i\overset{(\ref{opps100})}{\not\in} (\DI_{\delta}\cup \OI_{\delta})$ (recall definition (\ref{xcx10})) and so $\det(\na u(x_i))\geq \delta^{\frac{1}{100}}$ and hence 
\begin{eqnarray}
\label{band071}
\lt|u\lt(B_{\eta}(x_i)\backslash \BI_{\ep}\rt)\rt|
&\overset{(\ref{klk2}),(\ref{xcxxz101})}{\geq}&c\delta^{\frac{1}{100}}\eta^n -c\ep\eta^n  \nn\\
&\geq&c\delta^{\frac{1}{100}}\eta^n .
\end{eqnarray}
Now by (\ref{band02}), (\ref{kjk3}) (since $x_i\not\in \HI_{\rho_0}$ and $\eta\leq \frac{\rho_0}{2}$) we have 
\begin{equation}
\label{band077}
B_{8^{-1}\delta^{\frac{1}{100}}\eta}
(u(x_i))\subset u\lt( B_{\frac{\eta}{4}}(x_i)\rt). 
\end{equation}
So define 
\begin{eqnarray}
\label{band0124}
&~&A:=B_{\frac{\delta^{\frac{1}{100}}\eta}{64}}\lt(u(x_i)+e_1 \frac{\delta^{\frac{1}{100}} \eta}{16}\rt)\backslash u\lt(\BI_{\ep}\cap B_{\eta}(x_i)\rt)\nn\\
&~&\qd\qd\qd\text{ and }B:=B_{\frac{\delta^{\frac{1}{100}}\eta}{64}}\lt(u(x_i)-e_1 \frac{\delta^{\frac{1}{100}}\eta}{16}\rt)\backslash u\lt(\BI_{\ep}\cap B_{\eta }(x_i)\rt).
\end{eqnarray}
By (\ref{band0124}), (\ref{band077}) and the fact $u$ is injective on $B_{\eta}(x_i)$ (recall (\ref{ghf400.6}), (\ref{opps100}))
\begin{equation}
\label{kcka2}
A\cup B\subset u\lt(B_{\eta}(x_i)\backslash \BI_{\ep}\rt)
\end{equation}
and note 
\begin{equation}
\label{xcxzz1}
\mathrm{dist}(A,B)>\frac{\delta^{\frac{1}{100}}\eta}{64}.
\end{equation}
Now 
\begin{eqnarray}
\label{opps920}
\lt|A\rt|&\overset{(\ref{xcxxz101}),(\ref{band0124})}{\geq}& \frac{\delta^{\frac{n}{100}} \eta^n }{64^n}
-c\delta^{\frac{24n}{100}}\eta^n \nn\\
&\geq& c\delta^{\frac{n}{100}}\eta^n. 
\end{eqnarray}
In exactly the same way $\lt|B\rt|\geq c\delta^{\frac{n}{100}}\eta^n$. Now note 
\begin{eqnarray}
\label{band0.014}
&~&\eta^{-n}\ep^{\frac{n}{100}}\int_{u(B_{\eta}(x_i)\backslash \BI_{\ep})} \lt|l_R(z)-l_{R_i}(z)\rt| dz\nn\\
&~&\qd\qd\qd\qd\qd\qd \overset{(\ref{xcx11})}{\leq}
\eta^{-n}\int_{u(B_{\eta}(x_i)\backslash \BI_{\ep})} \lt|l_R(z)-l_{R_i}(z)\rt|\det(\na u(u^{-1}(z)))^{-1} dz\nn\\
&~&\qd\qd\qd\qd\qd\qd \overset{(\ref{cvc31})}{\leq} c \eta \ep^{2n}.
\end{eqnarray}
Let 
\begin{equation}
\label{cvcbb1}
U_A:=\lt\{z\in A:\lt|l_R(z)-l_{R^k_i}(z)\rt|>\eta\ep^{\frac{n}{2}}\rt\}.
\end{equation}
Notice 
\begin{eqnarray}
\eta \ep^{\frac{n}{2}}\lt|U_A\rt|&\leq&
\int_{A} \lt|l_R(z)-l_{R^k_i}(z)\rt| dz\nn\\
&\overset{(\ref{band0.014})}{\leq}& c\eta^{n+1} \ep^{\frac{199n}{100}} .\nn
\end{eqnarray}
So $\lt|U_A\rt|\leq c\eta^{n}\ep$, since $\ep<<\delta$, 
from (\ref{opps920}) $\lt|A\backslash U_A\rt|>0$ and we can pick $x_A\in A\backslash U_A$. In exactly 
the same way $x_B\in B\backslash U_B$. So $\lt|l_R(x_A)-l_{R^k_i}(x_A)\rt|\leq \eta\ep^{\frac{n}{2}}$ and $\lt|l_R(x_B)-l_{R^k_i}(x_B)\rt|\leq c\eta\ep^{\frac{n}{2}}$. Now 
$l_R(z)=R z+\alpha_R$ and $l_{R^k_i}(z)=R^k_i z+\alpha_{R_i}$ for some $\alpha_R,\alpha_{R_i}\in \R^n$, 
so  
$$
\lt|(R-R^k_i)x_A+(\alpha_R-\alpha_{R^k_i})\rt|\leq c\eta\ep^{\frac{n}{2}}
$$ 
and 
$$
\lt|(R-R^k_i)x_B+(\alpha_R-\alpha_{R^k_i})\rt|\leq c\eta\ep^{\frac{n}{2}}.
$$
Now taking one away from another $\lt|(R-R^k_i)(x_A-x_B)\rt|\leq c\eta\ep^{\frac{n}{2}}$. Note $\lt|x_A-x_B\rt|\overset{(\ref{xcxzz1})}{\geq} \frac{\delta^{\frac{1}{100}}}{64}\eta$ so 
\begin{equation}
\lt|(R-R^k_i)\frac{(x_A-x_B)}{\lt|x_A-x_B\rt|}\rt|\leq c\delta^{-\frac{1}{100}}\ep^{\frac{n}{2}}\leq c\ep^{\frac{n}{4}}\nn
\end{equation}
Therefor $\lt|R-R^k_i\rt|\leq c\ep^{\frac{n}{4}}$, this completes the proof of Step 2. \nl

\em Final step of Proof of Theorem \ref{TT2}.\rm 

Let $k>M_1$. For any $i\in\lt\{1,2,\dots m\rt\}\backslash (B^k_0\cup B^k_1\cup B_2\cup B_3)$
\begin{eqnarray}
\label{opps180}
\Xint{-}_{B_{\frac{\eta}{2}}(x_i)} \lt|\na v_k- R \na u_k\rt| dz&\leq&  \Xint{-}_{B_{\frac{\eta}{2}}(x_i)} \lt|\na v_k- R^k_i \na u_k\rt|
+ \lt|R^k_i \na u_k-R\na u_k\rt| dz\nn\\
&\overset{(\ref{cvc30.7})}{\leq}& c \ep^{2n}+\lt|R^k_i-R\rt|\Xint{-}_{B_{\frac{\eta}{2}}(x_i)} \lt|\na u_k\rt| dz\nn\\
&\overset{(\ref{vxx1})}{\leq}& c\ep^{2n}+c\ep^{\frac{n}{4}}\Xint{-}_{B_{\frac{\eta}{2}}(x_i)} \lt|\na u_k\rt| dz\nn\\
&\overset{(\ref{vvv7})}{\leq}& c\delta^{-1}\ep^{\frac{n}{4}}\leq c\ep^{\frac{n}{8}}.
\end{eqnarray}
Let 
\begin{equation}
\label{zband022.7}
\Pi':=\Pi\cap \lt(\bigcup_{\lt\{1,2,\dots m\rt\}\backslash (B^k_0\cup B^k_1\cup B_2\cup B_3)} 
B_{\frac{\eta}{2}}(x_i)\rt)
\end{equation} 
so 
\begin{eqnarray}
\label{opps101}
\lt|\Pi\backslash \Pi'\rt|&\overset{(\ref{zband022.7})}{\leq}& c\eta^n  \ca{B^k_0\cup B^k_1\cup B_2\cup B_3}\nn\\
&\overset{(\ref{cvc30})}{\leq}& c\sqrt{\gamma}+c\delta^{\frac{n}{2}}+c\lt|\DI_{\delta}\rt|.
\end{eqnarray}
Now recall from definition of $\Pi$ (\ref{opps100}) we have that 
\begin{eqnarray}
\label{opps106}
\lt|\Omega\backslash \Pi\rt|&\leq& \lt|\OI_{\delta}\rt|+\lt|\BI_{\delta}\rt|+\lt|\DI_{\delta}\rt|
+\lt|\EI_{\ep}\rt|+\lt|\Theta_{\rho_1}\rt|+\lt|\HI_{\rho_0}\rt|+\lt|N_{\delta}(\partial \Omega)\rt|\nn\\
&\overset{(\ref{opps101.5}),(\ref{ghf400.6}),(\ref{band0.018}),(\ref{kjk4})}{\leq}&c\delta^{\frac{1}{100}}+\lt|\BI_{\delta}\rt|+\lt|\DI_{\delta}\rt|+c\ep.
\end{eqnarray}
So putting (\ref{opps101}), (\ref{opps106}) together we have 
\begin{equation}
\label{opps600}
\lt|\Omega\backslash \Pi'\rt|\leq c\delta^{\frac{1}{100}}+\lt|\BI_{\delta}\rt|+\lt|\DI_{\delta}\rt|+c\ep+c\sqrt{\gamma}.
\end{equation}
So from (\ref{opps180}) and definition (\ref{zband022.7}) we have 
\begin{eqnarray}
\label{opps500}
\int_{\Pi'} \lt|\na v_k-R \na u_k \rt| dz&\leq& \sum_{i\in  \lt\{1,2,\dots m\rt\}\backslash (B^k_0\cup B^k_1\cup B_2\cup B_3)} 
\int_{B_{\eta}(x_i)} \lt|\na v_k-R \na u_k\rt| dz\nn\\
&\leq& c\ep^{\frac{n}{8}}.
\end{eqnarray}
To simplify notation let $\varsigma_k=\int_{\Omega} \lt|S(\na v_k)-S(\na u_k)\rt| dz$
\begin{eqnarray}
\label{opps201}
\lt|\int_{\Omega} \lt|\na v_k\rt| dz-\int_{\Omega} \lt|\na u_k\rt| dz\rt|&\leq&\int_{\Omega} \lt|\lt|\na v_k\rt|-\lt|\na u_k\rt|\rt| dz\nn\\
&\leq& c\int_{\Omega} \lt|\lt|S(\na v_k)\rt|-\lt|S(\na u_k)\rt|\rt| dz\nn\\
&\leq& c\int_{\Omega} \lt|S(\na v_k)-S(\na u_k)\rt| dz\nn\\
&\leq& c\varsigma_k.
\end{eqnarray}
Note also that 
\begin{eqnarray}
\label{opps200}
\int_{\Omega\backslash \Pi'} \lt|\na u_k\rt| dz&\leq& \lt(\int_{\Omega} \lt|\na u_k\rt|^n dz\rt)^{\frac{1}{n}}\lt|\Omega\backslash \Pi'\rt|^{\frac{n-1}{n}}\nn\\
&\overset{(\ref{opps600})}{\leq}& c\lt(\delta^{\frac{1}{100}}+\lt|\BI_{\delta}\rt|+\lt|\DI_{\delta}\rt|+\sqrt{\gamma}+\ep\rt)^{\frac{n-1}{n}}.
\end{eqnarray}
And 
\begin{eqnarray}
\label{opps300}
\int_{\Omega\backslash \Pi'} \lt|\na v_k\rt| dz&\overset{(\ref{opps201})}{\leq}&\int_{\Omega\backslash \Pi'} \lt|\na u_k\rt| dz+c\varsigma_k\nn\\
&\overset{(\ref{opps200})}{\leq}&c\lt(\delta^{\frac{1}{100}}+\lt|\BI_{\delta}\rt|+\lt|\DI_{\delta}\rt|+\sqrt{\gamma}+\ep\rt)^{\frac{n-1}{n}}+c\varsigma_k.
\end{eqnarray}
Thus 
\begin{eqnarray}
\label{opps350}
\int_{\Omega\backslash \Pi'} \lt|\na v_k-R \na u_k\rt| dz&\overset{(\ref{opps200}),(\ref{opps300})}{\leq}&c\lt(\delta^{\frac{1}{100}}+\lt|\BI_{\delta}\rt|+\lt|\DI_{\delta}\rt|+\sqrt{\gamma}+\ep\rt)^{\frac{n-1}{n}}+c\varsigma_k.
\end{eqnarray}
Putting this together with (\ref{opps500}) we have 
\begin{eqnarray}
&~&\int_{\Omega} \lt|\na v_k-R \na u_k\rt| dz \nn\\
&~&\qd\qd\qd\qd\leq c\lt(\delta^{\frac{1}{100}}+\lt|\BI_{\delta}\rt|+\lt|\DI_{\delta}\rt|+\sqrt{\gamma}+\ep^{\frac{1}{8}}\rt)^{\frac{n-1}{n}}+c\varsigma_k\text{ for all }k>M_1.\nn
\end{eqnarray}
Now recall $\ep<<\delta$ are $\gamma<<\ep$ and $\delta$ was chosen arbitrarily. So we have established (\ref{opps900}). $\Box$

%
%
%

\section{Counter example}
\label{counter}

\em Example 1. \rm  Let $Q_1:=\lt\{z:\lt|z\rt|_{\infty}<1\rt\}$. Define 
\begin{equation}
u(x_1,x_2,\dots x_n):=\lt\{\begin{array}{ll} (x_1,x_2 x_1, x_3, \dots x_n)&\text{ for } x_1>0\nn\\
(x_1,-x_2 x_1, x_3, \dots x_n)&\text{ for } x_1\leq 0\end{array}\rt.
\end{equation}
and for some $\theta\in (0,2\pi)$
\begin{equation}
v(x_1,x_2,\dots x_n):=\lt\{\begin{array}{ll} (x_1\cos \theta-x_1 x_2\sin \theta,x_1\sin\theta+x_1x_2\cos\theta,x_3, \dots x_n)&\text{ for } x_1>0\nn\\
(x_1,-x_2 x_1, x_3, \dots x_n)&\text{ for } x_1\leq 0\end{array}\rt.
\end{equation}
Note that for $x_1\leq 0$
$$
\na u(x)=\lt(\begin{matrix} 1&0&0&\dots 0\nn\\
-x_2&-x_1&0&\dots 0\nn\\
0&0&1&\dots 0\nn\\
\dots\nn\\
0&0&0&\dots 1\end{matrix}\rt)
$$
And for $x_1>0$
\begin{eqnarray}
\na v(x)&=&\lt(\begin{matrix} \cos\theta-x_2\sin\theta&-x_1\sin\theta&0&\dots 0\nn\\
\sin\theta+x_2\cos\theta&x_1\cos\theta&0&\dots 0\nn\\
0&0&1&\dots 0\nn\\
\dots\nn\\
0&0&0&\dots 1\end{matrix}\rt)\nn\\
&=&\lt(\begin{matrix} \cos\theta&-\sin\theta&0&\dots 0\nn\\
\sin\theta&\cos\theta&0&\dots 0\nn\\
0&0&1&\dots 0\nn\\
\dots\nn\\
0&0&0&\dots 1\end{matrix}\rt)\lt(\begin{matrix} 1&0&0&\dots 0\nn\\
x_2&x_1&0&\dots 0\nn\\
0&0&1&\dots 0\nn\\
\dots\nn\\
0&0&0&\dots 1\end{matrix}\rt).
\end{eqnarray}
Since $\na u(x)=\na v(x)$ for $x_1\leq 0$ it is clear there is no $R$ such that $\na v(x)=R \na u(x)$ for $x\in Q_1$. Now note that $\det(\na u(x))=x_1$ for all $x\in Q$ and $\lt|\na u(x)\rt|^n
=\lt((n-1)+x_2^2+x_1^2\rt)^{\frac{n}{2}}$ so defining $K(x):=\lt|\na u(x)\rt|^n/\det(\na u(x))=x_1^{-1}\lt((n-1)+x_2^2+x_1^2\rt)^{\frac{n}{2}}$. So it is clear that $\int_{Q_1} K(z) dz=\infty$ and 
thus it follows that Theorems \ref{T1} and \ref{TT2} are optimal for $n=2$.

\section{On the question of Sharpness of Theorem \ref{T1} and Theorem \ref{TT2}}
\label{nonocounter}

As mentioned the only known way of constructing a counter examples to Theorems \ref{T1} and \ref{TT2} is to 
take a function that squeezes down a domain into a shape whose interior consists of two disjoint pieces. In 
three dimensions in analogy with Example 1 of Section \ref{counter} we could consider squeezing the center of a cube 
to a line, in effect doing the squeezing only in the $x$ and $z$ variables. However in this case 
the calculations reduce to those of the two dimensional situation and it can be shown that for a wide class 
of mappings, squeezing down the center to a line implies that the mapping fails to 
have $L^1$ integrable dilatation. 

A more promising approach might be to consider mappings that squeeze down the center of a cylinder to a 
point. However Proposition \ref{nocounter} below will show, such examples (if they exist) can not be easily 
constructed. 

Let $R_{\theta}=\lt(\begin{matrix} \cos\theta & -\sin\theta  & 0 \\ \sin\theta & \cos\theta & 0 \\ 0 & 0 & 1\end{matrix}\rt)$ be a rotation 
around the $z$-axis. We say $\Omega$ is axially symmetric if $R_{\theta}\Omega=\Omega$ for every $\theta$. Now given a 
function $f:\Omega\rightarrow \R^3$  we say the function $f$ is axially symmetric if any axially symmetric subset $S\subset \Omega$ we have that $f(S)$ is axially symmetric. 

With a view to attempting to show sharpness of Theorems \ref{T1} and \ref{TT2} we would like to try and 
construct a function that squeezes $B_1(0)\times \lt[0,1\rt]$ in the center 
down to a point and use this to create a counter example to Theorems \ref{T1} and \ref{TT2} for functions whose dilatation are not 
$L^{p}$ for $p\geq n-1$. We say a function $g:\Omega\rightarrow \R$ (where $\Omega$ is axially symmetric) is a 
cylindrical product function if $g(r\cos \theta,r \sin\theta,z)=p_1(r)p_2(\theta)p_3(z)$ for functions $p_1,p_2,p_3$. A 
function $f:\Omega\rightarrow \R^3$ is a cylindrical product function if for each co-ordinate function is a cylindrical product function.

We will show that any axially symmetric orientation preserving cylindrical product function (whose coordinates 
satisfy certain monotonicity or convexity properties) that squeezes the cylinder down to a point does not have $L^1$ integrable dilatation.

\begin{a5}
\label{nocounter}
Let $f:W^{1,1}(B_1(0)\times \lt[0,1\rt]:\R^3)$ be a radially symmetric orientation preserving  cylindrical product function, i.e. 
there exists functions $w,v,g,h,l$ such that 
\begin{equation}
\label{aaeq2}
f(r\cos \theta,r \sin \theta, z)
=(w(z)v(r)\cos(g(\theta)),w(z)v(r)\sin(g(\theta)),h(z)l(r))
\end{equation}
for some functions $w,v,g,h,l$. Assume each of these functions $w,v,h,l$ are monotonic non-decreasing or non-increasing, 
$g$ is non decreasing and $w, h$ are either concave or convex.  If $f(B_1(0)\times \lt\{0\rt\})$ consists of a single point then 
\begin{equation}
\label{noeq1}
\int_{B_1(0)\times \lt[0,1\rt]} \frac{\|\na f\|^3}{\det(\na f)} dz=\infty.
\end{equation} 
\end{a5}
\em Proof of Proposition \ref{nocounter}. \rm  Suppose the proposition is false. So there exists a function $f$ satisfying the 
hypotheses and 
\begin{equation}
\label{eqeq1}
\int_{B_1(0)\times \lt[0,1\rt]} \frac{\|\na f\|^3}{\det(\na f)} dz<\infty.
\end{equation}
%
%

Let $u(\theta,r,z)=f(r\cos \theta,r \sin \theta, z)$, so 
$u(\theta,r,z)=(w(z)v(r)\cos(g(\theta)),w(z)v(r)\sin(g(\theta)),h(z)l(r))$. Since $w(0)=0$, this function is non 
decreasing, 
\begin{equation}
\label{zband036}
\frac{\partial w}{\partial z}(z)\geq 0\text{ for }a.e.\ z.
\end{equation}

Now we claim 
\begin{equation}
\label{zband034}
\frac{\partial v}{\partial r}(r)\geq 0\text{ for }a.e.\ r.
\end{equation}

So see this we argue as follows. Define 
$F:B_1(0)\times \lt[0,1\rt]\rightarrow \lt[0,1\rt]\times \lt[0,1\rt]$  by 
$F(x,y,z)=(\sqrt{x^2+y^2},z)$. Since $JAC(\na F)=\det(\na F \na F^T)=1$. Note by the Co-area formula 
$$
\infty>\CI=\int_{B_1(0)\times \lt[0,1\rt]} \lt|\na f\rt|^3 JAC(\na F) dz=\int_{\lt[0,1\rt]\times \lt[0,1\rt]}
\int_{F^{-1}(r,z)} \lt|\na f\rt|^3 dH^1 dr dz.
$$
Let $\delta\in (0,1)$ be some small number we decide on later. We can find a set $\GI\subset \lt[0,\delta\rt]\times \lt[0,1\rt]$ 
with $\lt|\GI\rt|\geq \frac{\delta}{2}$ and for any $(r,z)\in \GI$, $\int_{F^{-1}(r,z)} \lt|\na f\rt|^3 dH^1 \leq c\delta^{-1}$. 
Pick $(r,z)\in \GI$, by Holder's inequality we have 
\begin{equation}
\label{eqeq22}
\int_{F^{-1}(r,z)} \lt| \na f\rt| dH^1 \leq c\lt(\int_{F^{-1}(r,z)} \lt|\na f\rt|^3 dH^1 \rt)^{\frac{1}{3}}\delta^{\frac{2}{3}}\leq 
c\delta^{\frac{1}{3}}.
\end{equation}
How ever if $v$ is non increasing then even for very small $r$ we know $f(F^{-1}(r,z))$ must be the boundary of a 
a disc with radius $o(1)$, so we must have $H^1(F^{-1}(r,z))\sim o(1)$ which contradicts 
(\ref{eqeq22}). This establishes (\ref{zband034}).

Now we claim 
\begin{equation}
\label{zband030}
\frac{\partial h}{\partial z}(z)\geq 0\text{ for }a.e.\ z.
\end{equation}
To see this first assume $l$ is non constant, the for $r_1\not =r_2$ we have that 
$l(r_1)\not=l(r_2)$. Since $f(B_1(0)\times \lt\{0\rt\})$ consists of a single point we 
must have $h(0)l(r_1)=h(0)l(r_2)$ so we must have $h(0)=0$ and hence (\ref{zband030}) is 
established.

On the other hand if $l$ is constant then as $f$ is orientation preserving and 
\begin{eqnarray}
\label{zband31}
\na u&:=&
\lt(\begin{matrix}  \frac{dv}{dr} w\cos\circ g &    -w v\sin\circ g \frac{dg}{d\theta}     & \frac{dw}{dz} v\cos\circ g \\
 \frac{dv}{dr}w\sin\circ g  &    w v \cos\circ g \frac{dg}{d\theta}     &  \frac{dw}{dz} v\sin\circ g\\
 \frac{dl}{dr}h  & 0                 & \frac{dh}{dz} l \end{matrix}\rt)
\end{eqnarray}
and as $\det(\na u)=2\frac{\partial h}{\partial z} l \frac{\partial v}{\partial r} w^2 v \frac{\partial g}{\partial \theta}>0$ and 
so by (\ref{zband034}), (\ref{zband030}) is established.  

Now from (\ref{zband31}) we know
\begin{eqnarray}
\label{exuw1}
\|\na u(\theta,r,z)\|_{\infty}^3&\geq&\lt|\na u(\theta,r,z)\rt|^3\nn\\
&\geq& c\max\lt\{\lt|w(z) v(r) \frac{dg}{d\theta}(\theta)\rt|^3,\lt|\frac{dv(r)}{dr} w(z)\rt|^3,\lt|\frac{dw(z)}{dz} v(r)\rt|^3\rt.\nn\\
&~&\qd\qd\qd\lt.,\lt|\frac{dh}{dz}(z)l(r)\rt|^3,\lt|\frac{dl}{dr}(r)h(z)\rt|^3\rt\}.
\end{eqnarray}
And as by (\ref{zband030}), (\ref{zband034}) and the fact that $g$ is 
non decreasing $\frac{dv}{dr}\geq 0$, $\frac{dh}{dz}\geq 0$ and $\frac{dg}{d\theta}\geq 0$
\begin{eqnarray}
\label{exuw3}
\det(\na u(\theta,r,z))&=& -\frac{dh}{dz}(z)l(r)\frac{dv}{dr}(r)\frac{dg}{d\theta}(\theta) v(r) w(z)^2
+\frac{dw}{dz}(z)v(r)^2 w(z) \frac{dg}{d\theta}(\theta) \frac{dl}{dr}(r) h(z)\nn\\
&\overset{(\ref{zband036}),(\ref{zband030})}{\leq}& \frac{dw}{dz}(z)v(r)^2 w(z) \frac{dg}{d\theta}(\theta) \frac{dl}{dr}(r) h(z).
\end{eqnarray}
So by (\ref{exuw3}), (\ref{exuw1})
$$
\lt(\lt(\frac{dw}{dz}(z)\rt)^2/w(z)\rt) \lt| v(r)\rt|\lt(\lt| \frac{dg}{d\theta}(\theta) \frac{dl}{dr}(r) h(z)\rt|\rt)^{-1}
\leq \frac{\|\na u(\theta,r,z)\|^3}{\det(\na u(\theta,r,z))} 
$$
So 
\begin{eqnarray}
\label{exuw4}
\lt(\frac{dw}{dz}(z)\rt)^2/w(z)
\leq \lt|v(r)\rt|^{-1} \frac{\|\na u(\theta,r,z)\|^3}{\det(\na u(\theta,r,z))} \lt| \frac{dg}{d\theta}(\theta) \frac{dl}{dr}(r) h(z)\rt|.
\end{eqnarray}

\em Step 1. \rm We will show that there can not exists $\delta>0$ 
such that 
\begin{equation}
\label{eqeq2}
\sup\lt\{\lt|\frac{dw}{dz}(z)\rt|:z\in (0,\delta)\rt\}>\delta.
\end{equation}
\em Proof of Step 1. \rm By (\ref{exuw4}) we have that 
\begin{equation}
\label{exuw5}
\int_{\lt[0,\delta\rt]} w(z)^{-1} dz\leq c\int_{\lt[0,\delta\rt]}  \frac{\|\na u(\theta,r,z)\|^3}{\det(\na u(\theta,r,z))} dz 
\end{equation}
For $(\theta,r)\in \lt[0,2\pi\rt)\times \lt[0,1\rt]$ define $l_{(\theta,r)}:=\lt\{(\theta,r,z):z\in \lt[0,1\rt]\rt\}$. 
Let 
\begin{equation}
\label{eqeq3}
\GII:=\lt\{(\theta,r)\in \lt[0,2\pi\rt)\times \lt[0,1\rt]:\int_{l_{(\theta,r)}} \frac{\|\na u(\theta,r,z)\|^3}{\det(\na u(\theta,r,z))} dz<\infty\rt\}.
\end{equation}
Now by Fubini $\lt|\lt[0,2\pi\rt)\times \lt[0,1\rt]\backslash \GII\rt|=0$. So for 
any $(\theta,r)\in \GII$ by (\ref{exuw5}) 
we have that 
$$
\mu(A):=\int_{A} w(z)^{-1} dz
$$ 
forms a finite measure on the interval $\lt[0,\delta\rt]$. Let $\ep<<\delta$, 
now by Holder we have that 
\begin{eqnarray}
\log(w(\delta))-\log(w(\ep))&=&
\int_{\ep}^{\delta} \frac{d}{dz}(\log(w(z))) dz=\int_{\lt[0,\delta\rt]} \frac{dw}{dz}(z)/w(z) dz\nn\\
&=&\int_{\lt[\ep,\delta\rt]} \frac{dw}{dz}(z) d\mu (z)\leq \lt( \int_{\lt[\ep,\delta\rt]} \lt(\frac{dw}{dz}(z)\rt)^2 d\mu (z)  \rt)^{\frac{1}{2}}\nn\\
&=&   \lt( \int_{\lt[\ep,\delta\rt]} \lt(\frac{dw}{dz}(z)\rt)^2/w(z)  \rt)^{\frac{1}{2}} dz 
\overset{(\ref{exuw4}),(\ref{eqeq3})}{<} \infty.
\end{eqnarray}
Since $w(0)=0$ sending $\ep\rightarrow 0$ we have contradiction. \nl

\em Proof of the Proposition completed. \rm Firstly if $w$ is concave, we must have $\lim_{z\rightarrow 0} w'(z)>0$ since otherwise by the fact that 
$w(0)=0$ and $w$ is positive we would have $w\equiv 0$. So the existence of some $\delta$ satisfying (\ref{eqeq2}) follows and so we have a 
contradiction. If $w$ is convex and $\lim_{z\rightarrow 0} w'(z)>0$ then again its easy to see there exists $\delta$ satisfying (\ref{eqeq2}) 
and so have a contradiction. 

So the only remain case to consider the when $w$ is convex and $\lim_{z\rightarrow 0} w'(z)=0$. Pick $(\theta,r)\in \GI$, 
let 
$$
\AI:=\lt\{z\in \lt(0,1\rt):w(z)> h(z)\rt\}\text{ and }\BI:=\lt(0,1\rt)\backslash \AI.
$$
\begin{eqnarray}
\label{band0.91}
\frac{\|\na u(\theta,r,z)\|^3_{\infty}}{\det(\na u(\theta,r,z))}
&\overset{(\ref{exuw1}),(\ref{exuw3})}{\geq}& \lt|\frac{dw}{dz}(z)\rt|^2 v(r)/\lt(w(z) \frac{dg}{d\theta}(\theta) 
\frac{dl}{dr}(r)h(z)\rt)\nn\\
&\geq& c(\theta,r) \lt|\frac{dw}{dz}(z)\rt|^2/w(z)^2 \text{ for }z\in \AI
\end{eqnarray}
where $c(\theta,r):=v(r)/ (\frac{dg}{d\theta}(\theta) \frac{dl}{dr}(r))$.

Now by Young's inequality for $p=3, p'=\frac{3}{2}$ from (\ref{exuw1}) we have that 
\begin{eqnarray}
\| \na u(\theta,r,z)\|^3_{\infty}&\geq& \lt(l(r)^2 \lt(\frac{dh}{dz}(z)\rt)^2\rt)^{p'}+\lt(v(r)\frac{dw}{dz}(z)\rt)^p\nn\\
&\geq&c l(r)^2 v(r) \lt(\frac{dh}{dz}(z)\rt)^2\frac{dw}{dz}(z).\nn
\end{eqnarray}
So by (\ref{exuw3}) we have that 
\begin{eqnarray}
\label{eqeq16}
\frac{\|\na u(\theta,r,z)\|^3_{\infty}}{\det(\na u(\theta,r,z))}\geq d(\theta,r)\lt|\frac{dh}{dz}(z)\rt|^2/h(z)^2 \text{ for }z\in \BI
\end{eqnarray}
where $d(\theta,r):=c l(r)^2/\lt( \frac{dg}{d\theta}(\theta)\frac{dl}{dr}(r)\rt)v(r)$.

Now $\AI$ is open so is the union of a countable collection of disjoint open intervals, denote them $I_1,I_2,\dots$. Thus 
$\AI=\bigcup_{k=1}^{\infty} I_k$. We assume they have been ordered to that $\sup I_k\leq \inf I_j$ for $k>j$. 
Define $\alpha_{2k}, \alpha_{2k+1}$ to be the endpoints of $I_k$ where $\alpha_{2k+1}\leq \alpha_{2k}$, i.e. 
$I_k=(\alpha_{2k+1},\alpha_{2k})$. Define $J_k:=(\alpha_{2k+2},\alpha_{2k+1})$, so $[0,1]=\bigcup_{k=1}^{\infty} \overline{I_k}\cup \overline{J_k}$.

Hence by Holder for any $(\theta,r)\in \GII$ we have 
$\int_{\lt[0,1\rt]} \sqrt{ \frac{\|\na u(\theta,r,z)\|^3_{\infty}}{\det(\na u(\theta,r,z))} } dz<\infty$. So by (\ref{eqeq3}) we 
have 
\begin{eqnarray}
\infty&>&\int_{\AI}  \sqrt{\frac{\|\na u(\theta,r,z)\|^3_{\infty}}{\det(\na u(\theta,r,z))}}\nn\\
&\overset{(\ref{band0.91})}{\geq}& c(\theta,r)\int_{\AI} \lt|\frac{dw}{dz}(z)\rt|/w(z) dz. \nn
\end{eqnarray}
Now as $\lim_{k\rightarrow \infty} \int_{\alpha_{2k}}^1 \frac{d}{dz}\lt(\log(w(z))\rt) dz=\lim_{k\rightarrow \infty}(\log(w(1))-\log(w(\alpha_{2k})))=\infty$ we must have 
\begin{equation}
\label{eqeq11}
\lim_{k\rightarrow \infty}\int_{\BI\cap \lt[\alpha_{2k},1\rt]} \frac{d}{dz}\lt(\log(w(z))\rt) dz=\infty
\end{equation}
But note that $w(\alpha_{2k})=h(\alpha_{2k})$ and  $w(\alpha_{2k+1})=h(\alpha_{2k+1})$ for every $k$. So 
\begin{eqnarray}
\int_{\BI\cap \lt[\alpha_{2k},1\rt]} \frac{d}{dz}\lt(\log(w(z))\rt) dz&=&\sum_{i=1}^{k} 
\int_{\lt[\alpha_{2i},\alpha_{2i-1}\rt]}   \frac{d}{dz}\lt(\log(w(z))\rt) dz\nn\\
&=&\int_{\BI\cap \lt[\alpha_{2k},1\rt]}   \frac{d}{dz}\lt(\log(h(z))\rt) dz\nn\\
&\overset{(\ref{eqeq16})}{\leq}&c\int_{\BI\cap \lt[\alpha_{2k},1\rt]} \sqrt{ \frac{\|\na u(\theta,r,z)\|^3_{\infty}}{\det(\na u(\theta,r,z))}} dz\nn\\
&<&\infty \nn
\end{eqnarray}
which contradicts (\ref{eqeq11}). $\Box$

\end{document}